\documentclass[leqno,draft,11pt]{article}
\usepackage{amssymb,amsmath,latexsym,theorem,a4wide}

\def\sameenum{}

\ifx\mydefs\undefined\else \fi
\let\mydefs\relax

\ifx\slide\undefined
  \allowdisplaybreaks[2]
\fi

\ifx\loadcyr\undefined\else
\input cyracc.def
\makeatletter
 at 1\@ptsize pt
 at 1\@ptsize pt
 at 1\@ptsize pt
\makeatother

\fi

\def\gobble#1{}
\def\fixsup#1#2{{#1\let\dp\gobble\mathstrut}^#2_}

\def\bme{\hskip.75em\relax}

\let\EQ\Leftrightarrow

\def\iff{\quad\text{iff}\quad}
\let\LOR\bigvee
\let\ET\bigwedge
\let\TO\Rightarrow
\def\?{\mathbin?}

\let\model\vDash

\newbox\circlebox
\setbox\circlebox\hbox{$\bigcirc$}
\def\circled#1{%
  \setbox0\hbox to\wd\circlebox{\hss$#1$\hss}\wd0=0pt
  \box0\copy\circlebox}

\let\fii\varphi
\let\tet\vartheta
\let\ep\varepsilon

\let\roo\varrho
\ifx\slide\empty
  \def\greek#1{$\expandafter\greeknum\csname c@#1\endcsname$}
\else
  \def\greek#1{$\mathop{\boldsymbol{\expandafter\greeknum\csname c@#1\endcsname}}$}
\fi
\makeatletter
\def\greeknum#1{\ifcase#1\or\alpha\or\beta\or\gamma\or\delta\or\ep
      \or\digamma\or\zeta\or\eta\or\tet\or\iota\else\@ctrerr\fi}
\makeatother

\def\p#1{\langle#1\rangle}

\def\lh#1{\lvert#1\rvert}

\let\bez\smallsetminus

\let\sset\subseteq
\let\nsset\nsubseteq

\let\Sset\supseteq

\let\onto\twoheadrightarrow

\def\bigcupd{\mathop{\dot\bigcup}}
\def\pw#1{\mathcal P(#1)}
\let\nul\varnothing

\def\two{\mathbf2}

\newcommand\rpair[3][3em]{\mathrel{%
   \begin{matrix}%
     \strut\smash{\xrightonto{\hbox to#1{\hss$#2$\hss}}}\\[-1.7ex]%
     \strut\smash{\xleftembed[\hbox to#1{\hss$#3$\hss}]{}}%
   \end{matrix}}}
\makeatletter
\newcommand\xrightonto[2][]{\ext@arrow 0359\rightontofill{#1}{#2}}
\newcommand\xleftembed[2][]{\ext@arrow 3095\leftembedfill{#1}{#2}}
\def\leftembedfill{\arrowfill@\leftarrow\relbar\hookleftnoarrow}
\def\rightontofill{\arrowfill@\relbar\relbar\onto}
\def\hookleftnoarrow{\DOTSB\relbar\joinrel\rhook}
\makeatother

\mathchardef\#="2023 

\ifx\busspf\undefined\else
\usepackage{bussproofs}
\EnableBpAbbreviations

\fi

\ifx\symlasy\undefined

  \def\centdot#1{{%
    \setbox0\hbox{$\mathop{#1}$}\dimen0 \ht0
    \setbox0\hbox{$#1$}\advance\dimen0 -\ht0
    \setbox2\hbox to\wd0{\hss$\mathop{\cdot}$\hss}\wd2=0pt
    \lower\dimen0\box2\box0 }}
\else
  
  \def\centdot#1{{%
     \setbox0\hbox{$#1$}%
     \raise0.206\ht0\hbox to\wd0{\hss$\cdot$\hss}%
     \kern-\wd0 \box0 }}
\fi

\let\sls|
\def\up{\mathord\uparrow}
\def\down{\mathord\downarrow}
\def\Up{{\setbox0\hbox{$\uparrow$}%
         \lower\dp0\hbox to\wd0{\hss\vrule width4pt height.4pt\hss}%
         \kern-\wd0\box0}}
\def\UP{{\setbox0\hbox{$\uparrow$}%
         \lower\dp0\hbox to\wd0{\hss\vrule width4pt height.4pt\hss}%
         \kern-\wd0\copy0\kern-\wd0\raise.35ex\box0}}
\def\Down{{\setbox0\hbox{$\downarrow$}%
         \raise\ht0\hbox to\wd0{\hss\vrule width4pt depth.4pt\hss}%
         \kern-\wd0\box0}}

\newif\ifnadm

\def\doadm{\mathrel{%
   \setbox0 \hbox{$\mathop\vdash$}\dimen0 \ht0
   \setbox0 \hbox{$\vdash$}\advance\dimen0 -\ht0
   \vrule width.8\fontdimen8 \textfont3 height\ht0 depth\dp0
   \mkern-1mu
   \lower\dimen0 \hbox{$\vcenter{%
      \ifnadm
        \setbox0 \hbox{$\scriptstyle\sim\mathstrut$}%
        \hbox{\hbox to\wd0{\hss$\scriptstyle/$\hss}\kern-\wd0 \box0 }%
      \else
        \hbox{$\scriptstyle\sim\mathstrut$}%
      \fi}$}}}

\def\nrstyle#1#2#3{%
  \setbox0\hbox{$#1\bigcirc$}%
  \vcenter{\hbox to\wd0{\hss$#2#3$\hss}}%
  \kern-\wd0\box0 }

\DeclareMathOperator\dom{dom}
\DeclareMathOperator\rng{rng}

\DeclareMathOperator\id{id}

\DeclareMathOperator\cls{cl}

\DeclareMathOperator\var{var}

\ifx\slide\undefined
  
\fi

\DeclareMathOperator\pmf{Pmf}
\DeclareMathOperator\wgt{Wgt}
\DeclareMathOperator\inv{Inv}
\DeclareMathOperator\pol{Pol}
\DeclareMathOperator\sym{Sym}
\DeclareMathOperator\oker{oker}
\DeclareMathOperator\ocon{OCon}
\DeclareMathOperator\cst{c}

\def\st{\expandafter\hat}

\def\N{\mathbb N}
\def\Q{\mathbb Q}
\def\Z{\mathbb Z}

\def\GF#1{\mathbb F_{#1}}

\let\pres\triangleright
\let\npres\ntriangleright
\let\strc\mathbf

\mathcode`\*="0203
\let\ob\overline

\mathchardef\mhyphen="2D
\def\cput(#1)#2{\put(#1){\hbox to0pt{\hss#2\hss}}}
\def\txto{${}\to{}$}

\ifx\slide\undefined

\def\noproof{\leavevmode\unskip\bme\vadjust{}\nobreak\hfill$\qed$\par}
\let\qed\Box
\newenvironment{Pf}[1][]
  {\par\noindent\textit{Proof\optpar{#1}:}\bme\ignorespaces}
  {\noproof\pagebreak[2]\vskip\medskipamount\ignorespacesafterend}
\def\optpar#1{\ifx\relax#1\relax\else\ #1\fi}
\def\qedhere{\relax\ifmmode\eqno\qed\expandafter\aftergroup
                   \else\noproof\fi\noqed}
\def\noqed{\let\noproof\relax}

\theoremstyle{plain}
\ifx\shortthm\undefined
\newtheorem{Thm}{Theorem}[section]
\else
\newtheorem{Thm}{Theorem}
\fi
\newtheorem{Prop}[Thm]{Proposition}
\newtheorem{Cor}[Thm]{Corollary}
\newtheorem{Lem}[Thm]{Lemma}

\newtheorem{Prob}[Thm]{Problem}

\ifx\shortthm\undefined
\def\theCl{\arabic{Cl}}
\fi

\theorembodyfont\upshape
\newtheorem{Def}[Thm]{Definition}
\newtheorem{Rem}[Thm]{Remark}
\newtheorem{Exm}[Thm]{Example}

\newenvironment{Pf*}{\let\qed\qedCl\Pf}\endPf

\fi

\ifx\slide\undefined
  \usepackage[reftex]{theoremref}
\fi
\makeatletter
\def\thmref@flush{%
   \ifx\thmref@last\empty\else
      \ifthmref@comma, \thmref@finaltrue\fi \thmref@commatrue
      \thmref@last \ifx\thmref@stack\empty\else s\fi \thmref@num 0
      \let\do\thmref@one \thmref@stack
      \ifcase\thmref@num\or\space and\else\thmref@finaltrue, and\fi
      ~\ref{\thmref@head}\let\thmref@stack\empty\fi}
\def\thmref@one#1{\ifnum\thmref@num>0,\fi
   \space\ref{#1}\advance\thmref@num 1\relax}
\makeatother

\newif\iflinenumbers
\linenumberstrue

{\catcode`\^^I=13 \catcode`\^^M=13
\gdef\doalgo#1#2\end#{\hbox to\hsize{\hss \let^^I\qquad%
  \def\\^^M{\nobreak\hfil\break\vadjust{}\qquad}%
  \fboxsep1em \linenum0 %
  \fbox{\hsize#1\vbox{%
  \everypar{\advance\linenum1 %
      \hbox to1.2em{%
           \hss\iflinenumbers$\scriptstyle\the\linenum$\hskip.6em\fi}}%
  #2}}\hss}\end}}
\newcount\linenum

\def\key{\relax\ifmmode\expandafter\mathbf\else\expandafter\textbf\fi}

\def\allowhyphens{\nobreak\hskip0pt\relax}

\DeclareRobustCommand*\magiclparen{\ifmmode(\else\textup(\allowhyphens\fi}
\DeclareRobustCommand*\magicrparen{\ifmmode)\else\textup)\fi}
\let\lparen=(  \let\rparen=)
\def\magicparon{\catcode`\(\active\catcode`\)\active}
\def\magicparoff{\catcode`\(12 \catcode`\)12 }
\AtBeginDocument{\ifx\ifPreview\iftrue\else\magicparon\fi}
\magicparon
\let (=\magiclparen  \let )=\magicrparen

\ifx\sameenum\undefined
  \def\theenumi{\roman{enumi}}
  \ifx\enumup\undefined
    
  \else
    
  \fi

\else
  \def\theenumi{(\roman{enumi})}

\fi

\magicparoff

\mathchardef\comma=\mathcode`\,
{\catcode`\,=\active \gdef,{\comma\penalty\relpenalty}}

\ifx\slide\undefined

\providecommand\dedic{%
  \message{^^JWARNING: embed explicit dedication in the paper!^^J}%
  \thanks{Supported by
    grant IAA100190902 of GA AV \v CR, Center of Excellence CE-ITI under the grant
    P202/12/G061 of GA \v CR, and RVO: 67985840.}}

\author{Emil Je\v r\'abek\dedic\\[\medskipamount]
Institute of Mathematics of the Czech Academy of Sciences\\
\small \v Zitn\'a 25,
115\:67 Praha 1,
Czech Republic,
email: \texttt{jerabek@math.cas.cz}
}

\else\ifx

\author{Emil Je\v r\'abek}
\email{jerabek@math.cas.cz\\[-1em]http://math.cas.cz/\string~jerabek/}
\institution{Institute of Mathematics of the Czech Academy of Sciences, Prague}

\else 

\author[Emil Je\v r\'abek]{Emil Je\v r\'abek\\[\medskipamount]
   \scriptsize\texttt{jerabek@math.cas.cz}\\\texttt{http://math.cas.cz/\string~jerabek/}}
\institute{Institute of Mathematics of the Czech Academy of Sciences, Prague}

\fi\fi

\let\rng\relax
\DeclareMathOperator\rng{im}
\def\p#1{(#1)}

\title{Galois connection for multiple-output operations}

\def\dedic{\thanks{The research leading to these results has received funding from the European Research Council under
the European Union's Seventh Framework Programme (FP7/2007--2013)~/ ERC grant agreement no.~339691. The Institute of
Mathematics of the Czech Academy of Sciences is supported by RVO: 67985840.}}

\begin{document}
\maketitle
\begin{abstract}
It is a classical result from universal algebra that the notions of polymorphisms and invariants provide a Galois
connection between suitably closed classes (clones) of finitary operations $f\colon B^n\to B$, and classes (coclones)
of relations $r\sset B^k$. We will present a generalization of this duality to classes of (multi-valued, partial)
functions $f\colon B^n\to B^m$, employing invariants valued in partially ordered monoids instead of relations. In particular, our set-up
encompasses the case of permutations $f\colon B^n\to B^n$, motivated by problems in reversible computing.
\end{abstract}

\section{Introduction}

One of the pivotal notions in universal algebra is the concept of a \emph{clone}: a set of finitary operations $f\colon
B^n\to B$ on a base set~$B$, closed under composition (superposition), and containing all projections. A typical clone
is the set of term operations of an algebra with underlying set~$B$; in a sense, clones can be thought of as algebras
with their signature abstracted away. Apart from universal algebra, clones have numerous applications in logic and
computer science, many owing to the celebrated work of Post~\cite{post} who classified all clones on a two-element
base set.

An important tool in the study of clones is the \emph{clone--coclone duality}, originally discovered by
Geiger~\cite{geig} and Bodnarchuk, Kaluzhnin, Kotov, and Romov~\cite{bkkr:i,bkkr:ii}. In its most basic form valid for
finite base sets~$B$, it states that the natural preservation relation between operations $f\colon B^n\to B$ and
relations $r\sset B^k$ induces a Galois connection that provides a dual isomorphism of the lattice of clones to
the lattice of \emph{coclones}: sets of relations closed under definability by primitive positive (pp) formulas.

Many variants and extensions of the clone--coclone Galois connection appear in the literature; let us mention just a
few without attempting to make an exhaustive list. The use of relational invariants was pioneered by Krasner (see e.g.\
\cite{kras:gal-abs,kras:eng}), who developed a duality theory for endofunctions $f\colon B\to B$. The above-mentioned
seminal paper by Geiger~\cite{geig} considers not just the case of functions $f\colon B^n\to B$ for
finite~$B$, but also partial multi-valued functions, and he briefly indicates possible generalizations to infinite base
sets~$B$. Iskander~\cite{iskand} describes, for arbitrary~$B$, the closed classes dual to sets of partial operations
$B^n\to B$. Rosenberg~\cite{rosen:inf,rosen:part} develops, for arbitrary~$B$, a duality for clones of total or partial
operations $B^n\to B$ using infinitary relations as invariants. P\"oschel~\cite{posch} describes Galois-closed sets of
finitary relations and operations for arbitrary~$B$. Kerkhoff~\cite{kerkh} develops the duality in more general
categories than~Set. Couceiro~\cite{coucei} presents a duality for sets of heterogeneous (partial, multi-valued)
operations $A^n\to B$ with two base sets $A$ and~$B$.

In this paper, we present a new generalization of this Galois connection: rather than classes of operations
$f\colon B^n\to B$, we consider (partial, multi-valued) functions $f\colon B^n\to B^m$, where $m\ge0$; that is,
operations with multiple outputs, just like usual operations already may have multiple inputs. On the dual side, we
use functions $w\colon B^k\to\strc M$ valued in partially ordered monoids~$\strc M$ as invariants; note that
other variants of the Galois connection mentioned above generally stick to ordinary relations or something close in
spirit (e.g., infinitary relations, or pairs of relations).

Let us briefly explain the primary motivation for this generalization, which comes from the work of Aaronson, Grier,
and Schaeffer~\cite{aa-gr-sch} (anticipated in~\cite{aa:revclass}, where a preliminary version of some results of the
present paper were posted~\cite{ej:revclass}).

One way to model conventional computation is with Boolean circuits
$C\colon\{0,1\}^n\to\{0,1\}$. We may consider families of circuits using various sets of basic 
gates, and then the class of Boolean functions definable by circuits over some basis is a clone on~$B=\{0,1\}$. Such
clones were classified by Post~\cite{post}, and their description becomes relevant when discussing circuits for
restricted classes of functions.

In conventional computing, we may freely destroy or duplicate information: for example, on input $x,y$ we may
compute $x+y$. In \emph{reversible computing}~\cite{perum}, this is disallowed: computation is required to be (in principle)
step-by-step invertible. (The addition example above could be made reversible by computing the pair
$x,x+y$ instead.) One motivation to study reversible computing comes from consideration of physical constraints.
Since the underlying time-evolution operators of quantum mechanics are invertible, any physical realization of
irreversible computation must ``store'' the excess information in the form of side-effects; more specifically,
the second law of thermodynamics implies \emph{Landauer's principle}, which states that erasure of $n$~bits of
information incurs a certain increase of entropy proportional to~$n$ elsewhere in the system.
In practical terms, this means the computer must draw the corresponding amount of energy and dissipate it in the
environment as heat. In contrast, there are no known theoretical limits on the energy efficiency of reversible
computation.

Going a step further, models of \emph{quantum computing}~\cite{jaeger} are inherently reversible, as all computation
steps perform unitary (hence invertible) operators.

Now, reversible computation can be modelled using a suitable kind of
circuits made of reversible gates, which are permutations (bijections) $f\colon B^n\to B^n$. A class of functions
computable by circuits over some basis of reversible gates forms a ``clone'' of permutations satisfying
a handful of natural closure properties. The main result of~\cite{aa-gr-sch} is a description of all such
closed classes of permutations on $B=\{0,1\}$, which can be thought of as an analogue of Post's classification for reversible computing.

Thus, one of the design goals of the present paper is to develop a duality that applies as a special case to the kind
of closed classes of permutations $f\colon B^n\to B^n$ studied in~\cite{aa-gr-sch}, providing it with a broader
framework. Of course, it is highly desirable to include the classical case of functions $f\colon B^n\to B$ as well.
This naturally leads to consideration of functions $f\colon B^n\to B^m$ as a common generalization. We do not have a
particular reason to consider also partial multi-valued functions, except that it happens to work; in fact, the basic
Galois connection is easier to study for partial multi-valued functions, while requirements of totality bring in extra
complications.

The paper is organized as follows. We begin with a handful of motivating examples in Section~\ref{sec:examples}. We
recall some preliminary facts about partially ordered structures and Galois connections in
Section~\ref{sec:ordered-structures}. In Section~\ref{sec:clones-coclones}, we present our main Galois connection
between partial multi-valued multi-output functions and pomonoid-valued weight functions in its most general form. In
Section~\ref{sec:restricted-cases}, we discuss variants of our Galois connection for restricted classes of multi-output
functions or weights that might be important for applications, in particular, we give Galois connections for classes of
\emph{total} multi-output functions in Section~\ref{sec:totality-conditions}, and for classes of permutations closed
under the ancilla rule as in~\cite{aa-gr-sch} in Section~\ref{sec:ancillas}. We discuss the merits of using weights in
subdirectly irreducible pomonoids in Section~\ref{sec:subd-irred-weights}, including a brief description of finitely
generated subdirectly irreducible commutative monoids as classified by Grillet~\cite{grillet}. A few concluding words
are included in Section~\ref{sec:conclusion}.

\section{Initial examples}\label{sec:examples}

The Galois connection we are going to introduce involves more complicated invariants on the ``coclone'' side in
contrast to the classical clone--coclone duality and many of its known generalizations: rather than some form of
relations, we need to use functions valued in partially ordered monoids. Before we get to the formal business, we will
present a few examples to show that this is in fact a necessary move that follows the nature of multiple-output
functions, and to have something easily graspable in mind to illustrate the subsequent abstract definitions.

First, let us recall the preservation relation from the standard Galois connection. A function $f\colon B^n\to B$
\emph{preserves} a relation $r\sset B^k$ if for every pair of matrices\footnote{We index tuples and matrices starting
from~$0$, so that e.g., an $n$-tuple is written as $\p{x_0,\dots,x_{n-1}}$. Accordingly, we write $\N=\{0,1,2,\dots\}$,
and index variables such as $i$ and~$j$ are implicitly taken in~$\N$; for instance, a quantifier over $i<n$ stands for $i=0,\dots,n-1$.}
\[
a=\left(\begin{array}{llcl}
a^0_0&a^0_1&\cdots&a^0_{n-1}\\
a^1_0&a^1_1&\cdots&a^1_{n-1}\\
\:\vdots&\:\vdots&\ddots&\:\vdots\\
a^{k-1}_0&a^{k-1}_1&\cdots&a^{k-1}_{n-1}
\end{array}\right)
\in B^{k\times n},\qquad
b=\left(\begin{array}{l}b_0^0\\b_0^1\\\vdots\\b_0^{k-1\!\!\!}\end{array}\right)
\in B^{k\times 1},
\]
if all rows represent values of $f$:
\[f(a_0^j,\dots,a_{n-1}^j)=b_0^j,\quad j<k,\]
and if all columns of $a$ are in~$r$:
\[\p{a_i^0,\dots,a_i^{k-1}}\in r,\quad i<n,\]
then the (unique in this case) column of $b$ is in~$r$ as well:
\[\p{b_0^0,\dots,b_0^{k-1}}\in r.\]
In order to make the notation more concise, we will denote a matrix~$a$ as above by $(a_i^j)_{i<n}^{j<k}$; its rows and
columns will be denoted $a^j=\p{a_i^j:i<n}$ and $a_i=\p{a_i^j:j<k}$, 
respectively. Thus, we can state the definition of the preservation relation as follows: for all $a\in B^{k\times n}$
and $b\in B^{k\times 1}$, if $f(a^j)=b_0^j$ for every $j<k$, and $a_i\in r$ for every $i<n$, then $b_0\in r$. (So far
it would have been simpler to treat $b$ as a vector rather than a $1$-column matrix, but we will shortly realize this is an
artifact of $f$ having unary output.)

We will now have a look at some of the closed classes of reversible operations on two-element base set, considered
in~\cite{aa-gr-sch}.
\begin{Exm}\th\label{exm:clone-cons}
The class of all \emph{conservative} permutations $f\colon B^n\to B^n$ for $B=\{0,1\}$: that is, bijective functions such
that $f(x)$ and~$x$ have the same Hamming weight (= the number of $1$s) for any $x\in B^n$. We can express this
condition as follows: a permutation $f\colon B^n\to B^n$ is conservative iff
\begin{equation}\label{eq:12}
f(a^0)=b^0\implies\sum_{i<n}a^0_i=\sum_{i<n}b_i^0
\end{equation}
holds for all matrices $a,b\in B^{1\times n}$, where the sum is computed in the integers, viewing $B$ as a subset of~$\N$.
\end{Exm}
\begin{Exm}\th\label{exm:clone-modpres}
(Still $B=\{0,1\}$.) The class of all \emph{mod-$c$-preserving} permutations $f\colon B^n\to B^n$ for a constant $c>1$: that
is,
permutations such that the Hamming weights of $f(x)$ and~$x$ are congruent modulo~$c$ for any $x\in B^n$. A permutation
$f\colon B^n\to B^n$ is mod-$c$-preserving iff it satisfies the property~\eqref{eq:12} for all matrices $a,b\in
B^{1\times n}$, where the sum is now computed in the group $\Z/c\Z$.
\end{Exm}

The previous two examples show that multi-output functions, and specifically permutations, can ``count''---they are
capable of preserving numerical quantities associated with the input (as opposed to yes/no properties as given
by relations $r\sset B^k$), expressible as certain ``sums'' over the input elements.

In the general definition in Section~\ref{sec:clones-coclones}, we will actually employ multiplicative notation, so the
``sums'' will be written as ``products''; this is partly to emphasize that we will allow the aggregation operation to be
noncommutative, but mostly because we will at some point need to make this ``multiplication''
operation interact in a ring-like fashion with another kind of ``addition''. In any case, this is just a matter of
notation.

The next example also appears in~\cite{aa-gr-sch} when restricted to permutations, but we will state it for general
functions in order to showcase new features. In particular, in the first two examples, the ``sums'' (to become
``products'') were preserved by equality; but in the general case, they will only be preserved by inequality.
\begin{Exm}\th\label{exm:clone-aff}
($B=\{0,1\}$) 
The class of \emph{affine functions} $f\colon B^n\to B^m$, i.e., $f(x)=Ax+c$ for some $A\in B^{m\times n}$ and~$c\in
B^m$, identifying $B$ with the field~$\GF2$. In order to characterize this class in a similar spirit to
\th\ref{exm:clone-cons,exm:clone-modpres}, notice first that $f$ is affine if and only if each of its $m$~components
$f_i\colon\GF2^n\to\GF2$ is affine. This gets us in the realm of the classical clone--coclone duality: we know that
$f_i$ is affine iff it preserves the relation
\[r=\{\p{a^0,a^1,a^2,a^3}\in\GF2^4:a^0+a^1+a^2+a^3=0\}\]
(the sum being computed in~$\GF2$). Thus, let $w\colon\GF2^4\to\{0,1\}$ denote the characteristic function of~$r$, i.e.,
\[w(x^0,x^1,x^2,x^3)=x^0+x^1+x^2+x^3+1\]
(still using $\GF2$ addition). Then $f\colon B^n\to B^m$ is affine iff for all $a\in B^{4\times n}$ and $b\in
B^{4\times m}$,
\[\forall j<4\,f(a^j)=b^j\land\forall i<n\,w(a_i)=1\implies\forall i<m\,w(b_i)=1.\]
We can recast this as preservation of a product: $f$ is
affine iff for all $a\in B^{4\times n}$ and $b\in B^{4\times m}$,
\[\forall j<4\,f(a^j)=b^j\implies\prod_{i<n}w(a_i)\le\prod_{i<m}w(b_i).\]

In a similar way, any classical relational invariant $r\sset B^k$ can be expressed as preservation of a product of
certain functions valued in the two-element meet-semilattice $\p{\{0,1\},1,{\cdot},\le}$.
\end{Exm}

Invariants of a similar syntactic shape as above can be defined for functions $w\colon B^k\to\strc M$, where $\strc M$
is a structure in which we can compute products of finite sequences of elements, and we have a suitable order
relation. That is, we are led to the class of \emph{partially ordered monoids}.

\section{Preliminaries}\label{sec:ordered-structures}

\subsection{Ordered structures}

Since partially ordered monoids and other ordered structures are omnipresent in this paper, this section presents a
summary of some background information on such structures.

We might as well start from the beginning, even though the reader must have seen monoids and partial orders before. So,
recall that a binary relation ${\le}\sset X\times X$ is a \emph{preorder} on~$X$ if it is reflexive and transitive. If
it is additionally antisymmetric ($x\le y\land y\le x\to x=y$), it is a \emph{partial order}.

Let $\le$ be a partial order on~$X$. A set $Y\sset X$ is a \emph{down-set} if $x\le y$ and $y\in Y$ implies $x\in Y$,
and an \emph{up-set} if it satisfies the dual condition. For any $Y\sset X$, $Y\down$ denotes the generated down-set
$\{x\in X:\exists y\in Y\,x\le y\}$, and $Y\up$ the generated up-set $\{x\in X:\exists y\in Y\,y\le x\}$.

A structure $\strc M=\p{M,1,{\cdot}}$ in a signature with one constant, and one binary operation, is a
\emph{monoid} if $\cdot$ is associative, and $1$ is a two-sided unit ($1\cdot x=x\cdot 1=x$). We often write just $xy$
for $x\cdot y$. Associativity allows us to unambiguously refer to iterated products
\[x_0x_1\cdots x_{n-1}=\prod_{i<n}x_i.\]
This product is understood to be~$1$ if $n=0$. We will also sometimes write monoids in the additive signature
$\p{M,0,+}$, particularly when commutative.

A \emph{partially ordered monoid} (\emph{pomonoid} for short) is a structure $\strc M=\p{M,1,{\cdot},{\le}}$ such that
$\p{M,1,{\cdot}}$ is a monoid, and $\le$ is a partial order on~$M$ compatible with multiplication, i.e., satisfying
\begin{align*}
x\le y&\to xz\le yz,\\
x\le y&\to zx\le zy
\end{align*}
for all $x,y,z\in M$.

Apart from the class of pomonoids, we will also need to work e.g.\ with its subvarieties, and expansions such as
semirings. Thus, it will be helpful to have a general framework for ordered structures.
While we will assume the reader is familiar with basic notions from universal algebra---such as varieties, equational logic,
and subdirect irreducibility---in the standard set-up of purely algebraic structures (see e.g.~\cite{bur-san}), the
corresponding theory for partially ordered structures and inequational logic is much less commonly known, hence we
review the relevant concepts below. The results we need can be found in~\cite{pigo:povar}, but various parts of the
theory appear e.g.\ in \cite{bloom:povar,malc:alg-sys}.

Let us fix an algebraic signature~$\Sigma$ (in our application, it will most often be the signature of monoids with extra
constants). A \emph{partially ordered $\Sigma$-algebra} (short: poalgebra) is a $\Sigma$-algebra~$\strc A$ endowed with
a partial order~$\le_\strc A$ that makes
all functions from~$\Sigma$ monotone (nondecreasing) in every argument\footnote{More generally, we could specify for
each argument a \emph{polarity} indicating whether the function is nondecreasing or nonincreasing in the given
argument, see~\cite{pigo:povar}.}. Homomorphisms, subalgebras, products, and restricted products of poalgebras are
defined in the expected way, using the corresponding algebraic and order-theoretic notions on the respective parts of
the structure.

If $\phi\colon\strc A\to\strc B$ is a homomorphism between two poalgebras, its \emph{order kernel} is the relation
$\oker(\phi)=\{\p{a,a'}\in A^2:\phi(a)\le_\strc B\phi(a')\}$. The order kernel is an \emph{invariant preorder}
on~$\strc A$: a preorder $\preceq$ extending $\le_\strc A$ such that for each $f\in\Sigma$, $f^\strc A$ is
nondecreasing with respect to~$\preceq$ in every argument. Conversely, let $\preceq$ be an invariant preorder on~$\strc
A$. The relation ${\sim}={\preceq}\cap{\succeq}$ is a congruence of the algebraic part of~$\strc A$, hence we can form the
quotient structure $\strc B=\strc A/{\sim}$, and make it a poalgebra ordered by ${\preceq}/{\sim}$. The quotient map
$\phi\colon\strc A\to\strc B$ is a (surjective) homomorphism, and $\oker(\phi)={\preceq}$. We will denote $\strc B$ as
$\strc A/{\preceq}$.

Since invariant preorders are closed under arbitrary intersections, they form an algebraic complete lattice
$\ocon\strc A$, which
plays much the same role as the congruence lattice for unordered algebras.

A poalgebra $\strc A$ is a \emph{subdirect product} of a family of poalgebras $\{\strc A_i:i\in I\}$ if there exists an
embedding
\begin{equation}\label{eq:11}
\fii\colon\strc A\to\prod_{i\in I}\strc A_i
\end{equation}
such that $\pi_i\circ\fii\colon\strc A\to\strc A_i$ is surjective for each $i\in I$, where $\pi_i$ denotes
projection to the $i$th coordinate. A poalgebra $\strc A$ is \emph{subdirectly
irreducible} if for every subdirect product~\eqref{eq:11}, there exists $i\in I$ such that $\pi_i\circ\fii$ is an
isomorphism. An intrinsic characterization is that a poalgebra $\p{\strc A,\le_\strc A}$ is subdirectly irreducible iff it has a least invariant preorder properly
extending~$\le_\strc A$. Every poalgebra can be written as a subdirect product of subdirectly irreducible poalgebras.

A class of $\Sigma$-poalgebras is a \emph{(partially ordered) variety} if it can be axiomatized by a set of
(implicitly universally quantified) inequalities $t(\vec x)\le s(\vec x)$, where $t$ and~$s$ are $\Sigma$-terms, over
the theory of all $\Sigma$-poalgebras. Equivalently, a povariety is a class of poalgebras closed under subalgebras,
products, and homomorphic images (i.e., quotients by invariant preorders). If a poalgebra~$\strc A$ from a variety~$V$ is a
subdirect product of a family of poalgebras $\{\strc A_i:i\in I\}$, then each $\strc A_i$ is in~$V$, too; thus, every
$\strc A\in V$ is a subdirect product of subdirectly irreducible poalgebras from~$V$.

More generally, a class of poalgebras is a \emph{(partially ordered) quasivariety} if it can be axiomatized by a set of \emph{quasi-inequalities}
\[\ET_{i<n}t_i(\vec x)\le s_i(\vec x)\to t(\vec x)\le s(\vec x),\]
where $n\in\omega$, and $t,s,t_i,s_i$ are terms. Equivalently, a class of poalgebras is a quasivariety iff it is closed
under isomorphic images, subalgebras, products, and ultraproducts (or: restricted products).

If $Q$ is a quasivariety, and $\strc A$ a poalgebra, a \emph{$Q$-preorder} is an invariant preorder $\preceq$ on~$\strc
A$ such that $\strc A/{\preceq}\in Q$. The set of $Q$-preorders is closed under arbitrary intersections, hence it forms
an algebraic complete lattice~$\ocon_Q\strc A$. (If $Q$ is a povariety, then $\ocon_Q\strc A$ is a principal filter of
$\ocon\strc A$.) An algebra $\strc A\in Q$ is \emph{subdirectly irreducible relative to~$Q$} if for every subdirect
product~\eqref{eq:11} with $\{\strc A_i:i\in I\}\sset Q$, there is $i\in I$ such that $\pi_i\circ\fii$ is an
isomorphism. Equivalently, there is a least $Q$-preorder properly extending~$\le_\strc A$. Every $\strc A\in
Q$ can be written as a subdirect product of a family of poalgebras $\strc A_i\in Q$, subdirectly irreducible relative
to~$Q$. 

An unordered $\Sigma$-algebra $\strc A$ can be identified with the poalgebra $\p{\strc A,{=_\strc A}}$, which we will call
\emph{trivially ordered}. Notice that the class of trivially ordered poalgebras is a quasivariety, being axiomatized by
$x\le y\to y\le x$.

\subsection{Galois connections}\label{sec:galois-connections}

In order to ensure that we are all on the same page, let us also recall basic facts about closure operators and Galois
connections.

Let $\p{P,\le}$ be a partially ordered set. A \emph{closure operator} on~$P$ is a function $\cls\colon P\to P$ such 
that
\begin{equation}\label{eq:16}
X\le\cls Y\iff\cls X \le\cls Y
\end{equation}
for all $X,Y\in P$. (In our applications, $P$ will be typically the powerset lattice $\p{\pw
U,\sset}$ for some set~$U$; in this case, we will also call $\cls$ a closure operator \emph{on~$U$} for simplicity.) The
definition~\eqref{eq:16} is equivalent to the conjunction of the three conditions
\begin{gather*}
X\le\cls X,\\
X\le Y\implies\cls X\le\cls Y,\\
\cls X=\cls\cls X.
\end{gather*}
If $\cls$ is a closure operator, an $X\in P$ is \emph{closed} if $X=\cls X$. Let $C\sset P$ be the collection of all
closed elements. The definition of a closure operator implies that for any $X\in P$, $\cls X$ is the least closed element
above~$X$:
\begin{equation}\label{eq:13}
\cls X=\min\{Y\in C:X\le Y\}.
\end{equation}
That is, we can recover $\cls$ from~$C$.

A subset $C\sset P$ with the property that \eqref{eq:13} exists
for all $X\in P$ is called a \emph{closure system} on~$P$. For any closure system~$C$, the function $\cls$ defined
by~\eqref{eq:13} is a closure operator. The constructions of
closure systems from closure operators and vice versa described above are mutually inverse, hence the two definitions
can be considered different presentations of the same concept.
If $P$ is a complete lattice (e.g., a powerset lattice), closure systems have a simpler characterization: $C\sset P$ is
a closure system iff it is closed under arbitrary meets (including the empty meet, which yields the top element
of~$P$). In particular, any such closure system is itself a complete lattice.

Now, let $\p{P,\le}$ and $\p{Q,\preceq}$ be two partially ordered sets (again, typically powersets for us). A
\emph{Galois connection} between $P$ and~$Q$ is a pair of mappings $F\colon P\to Q$ and $G\colon Q\to P$ such that
\begin{equation}\label{eq:14}
X\le G(Y)\iff Y\preceq F(X)
\end{equation}
for all $X\in P$ and $Y\in Q$. This is equivalent to the conditions
\begin{gather*}
X\le G(F(X)),\\
Y\preceq F(G(Y)),\\
X\le X'\implies F(X')\preceq F(X),\\
Y\preceq Y'\implies G(Y')\le G(Y)
\end{gather*}
for $X,X'\in P$ and~$Y,Y'\in Q$.

A Galois connection as above induces closure operators $\cls_P=G\circ F\colon P\to P$ and $\cls_Q=F\circ G\colon Q\to
Q$ on $P$ and~$Q$, respectively. Elements of $P$ or~$Q$ are called \emph{closed} with respect to the Galois connection
(for short: \emph{Galois-closed}) if they are closed under $\cls_P$ or~$\cls_Q$, respectively.

The images of $F$ and~$G$ consist of closed elements. Moreover, $F$ and~$G$ restricted to the collections of
closed elements of $P$ and~$Q$ (respectively) are mutually inverse antitone isomorphisms.

The Galois connections discussed in this paper  arise by means of the following simple but powerful observation.
Let $U$ and~$V$ be sets, and $R\sset U\times V$ a binary relation. Then the mappings
\begin{equation}\label{eq:15}
\begin{aligned}
F(X)&=\{v\in V:\forall x\in X\,R(x,v)\},\\
G(Y)&=\{u\in U:\forall y\in Y\,R(u,y)\}
\end{aligned}
\end{equation}
form a Galois connection between the powerset lattices $\p{\pw U,\sset}$ and $\p{\pw V,\sset}$: indeed, either of
$X\sset G(Y)$ and $Y\sset F(X)$ is equivalent to the symmetric condition
\[\forall x\in X\,\forall y\in Y\,R(x,y),\]
hence \eqref{eq:14} holds.

We will not be too fussy about applying the concepts above to proper classes instead of sets, even though a
``powerclass lattice'' is not an honest object. In such cases, we will
take care to only work with maps from classes to classes that are definable, and refrain from problematic steps like
explicit quantification over subclasses, so it should be straightforward to formalize all our reasoning in ZFC.

\section{Multiple-output clones and coclones}\label{sec:clones-coclones}

Let us fix a base set $B$ for the rest of the paper. (While some results will only apply if $B$ is finite, the general
set-up allows arbitrary $B$. It may even be empty.)

As we already mentioned, we are going to study classes of \emph{partial multi-valued functions (pmf)} from
$B^n$ to~$B^m$ for some $n,m\in\omega$. Formally speaking, a pmf from $X$ to~$Y$ is just a relation $f\sset X\times Y$;
however, we view it as a nondeterministic operation that maps $x\in X$ to one of the values $y\in Y$ such that
$\p{x,y}\in f$ (if any). In order to stress this interpretation (while at the same time distinguishing it in notation from
``proper'' functions), we will use the symbol $f\colon X\TO Y$ to mean that $f$ is a pmf from $X$ to~$Y$, and we will
write $f(x)\approx y$ for $\p{x,y}\in f$.
\begin{Def}\th\label{def:pmf}
For any $n,m\in\omega$, let $\pmf_{n,m}$ denote the set of all partial multi-valued functions (pmf) from $B^n$ to~$B^m$
(that is, formally, $\pmf_{n,m}=\pw{B^n\times B^m}$). We also put $\pmf=\bigcupd_{n,m}\pmf_{n,m}$.
\end{Def}

Following the examples in Section~\ref{sec:examples}, we will characterize suitably closed classes of pmf by
preservation of certain ``weighted products'' in pomonoids. Thus, our invariants will be the following kind of objects:

\begin{Def}\th\label{def:wgt}
If $k\in\omega$, $\wgt_k$ denotes the class of $k$-ary \emph{weight functions,} i.e., mappings $w\colon B^k\to\strc M$
where $\strc M=\p{M,1,\cdot,\le}$ is any pomonoid. Let $\wgt=\bigcupd_k\wgt_k$.
\end{Def}

Without further ado, here is our fundamental preservation relation. (Keep in mind the notational conventions for matrices
from Section~\ref{sec:examples}.)
\begin{Def}\th\label{def:pres}
Let $n,m,k\in\omega$. A pmf $f\colon B^n\TO B^m$ \emph{preserves} a weight $w\colon B^k\to\strc M$, written $f\pres w$, if
the following holds for all $a=(a_i^j)_{i<n}^{j<k}\in B^{k\times n}$ and $b=(b_i^j)_{i<m}^{j<k}\in B^{k\times m}$:
\begin{equation}\label{eq:1}
\forall j<k\,f(a^j)\approx b^j\implies\prod_{i<n}w(a_i)\le\prod_{i<m}w(b_i).
\end{equation}
When $f\pres w$, we also say that $w$ is an \emph{invariant} of~$f$, and $f$ is a \emph{polymorphism} of~$w$. If
$C\sset\pmf$ and $D\sset\wgt$, we will write $C\pres D$ as a shorthand for $\forall f\in C\,\forall w\in D\,f\pres
w$.

Using~\eqref{eq:15}, the preservation relation induces a Galois connection between sets $C\sset\pmf$, and classes
$D\sset\wgt$:
\begin{align*}
\inv(C)&=\{w\in\wgt:C\pres w\},\\
\pol(D)&=\{f\in\pmf:f\pres D\}.
\end{align*}
\end{Def}
\begin{Cor}
$\pol\circ\inv$ and $\inv\circ\pol$ are closure operators on $\pmf$ and $\wgt$, respectively. $\inv$ and~$\pol$
are mutually inverse antitone isomorphisms between Galois-closed subsets of~$\pmf$, and Galois-closed subclasses
of~$\wgt$. $\inv(C)$ and~$\pmf(D)$ are Galois-closed for each $C\sset\pmf$ and~$D\sset\wgt$.
\noproof\end{Cor}

Our fundamental task in this section is to find an intrinsic characterization of Galois-closed subsets of~$\pmf$, and
subclasses of~$\wgt$.

Before we do that, we need to clarify one minor issue. The description of Galois-closed sets involves a
condition that says, roughly, that whether a pmf belongs in such a set depends only on its finite
parts. In the literature, this condition appears in several formulations under several names, such as \emph{local closure}. We
prefer to think of it as a topological closure property, however we include a few equivalent forms of the condition
below for the benefit of the reader.

\begin{Def}\th\label{def:two}
Let $\two$ be the set $\{0,1\}$; $\two_H$ denotes $\two$ endowed with the discrete Hausdorff topology, and $\two_S$
denotes $\two$ endowed with the Sierpi\'nski topology where $\{1\}$ open, but $\{0\}$ is not.
\end{Def}

\begin{Lem}\th\label{lem:sie}
Let $A$ be a family of subsets of~$X$, identified with their characteristic functions (i.e., elements of $\two^X$). The following are equivalent.
\begin{enumerate}
\item\label{item:14} $A$ is closed in~$\two_S^X$.
\item\label{item:15} $A$ is closed in~$\two_H^X$ and closed under subsets.
\item\label{item:16} $A$ is closed under directed unions and subsets.
\item\label{item:26} $A$ is of finite character: i.e., a $Y\sset X$ is in~$A$ iff all finite subsets of~$Y$ are in~$A$.
\end{enumerate}
\end{Lem}
\begin{Pf}

\ref{item:14}\txto\ref{item:15}: $A$ is closed in~$\two_H^X$ as $\two_H$ is
finer than~$\two_S$. Moreover, the basic open sets in~$\two_S^X$ are of the form $\{Y:Y\Sset Y_0\}$ where $Y_0\sset X$
is finite, hence every $\two_S^X$-open set is closed upwards, and every closed set is closed downwards.

\ref{item:15}\txto\ref{item:16}:
Let $S\sset A$ be directed, and $Y=\bigcup S$. For every finite $X_0\sset X$, there is $Y'\in S$ such that $Y'\Sset
Y\cap X_0$ by directedness, which implies $Y'\cap X_0=Y\cap X_0$. Thus, every $\two_H^X$-open neighbourhood of~$Y$
intersects~$A$, whence $Y\in A$.

\ref{item:16}\txto\ref{item:26}:
Left to right follows from closure under subsets. Right to left: $Y$ is a directed union of its finite subsets.

\ref{item:26}\txto\ref{item:14}:
Let $Y\notin A$. By finite character, there is a finite $Y_0\sset Y$ such that
$Y_0\notin A$, and then the basic $2_S^X$-open neighbourhood $\{Y':Y'\Sset Y_0\}$ of~$Y$ is disjoint from~$A$.
\end{Pf}

\begin{Rem}\th\label{rem:closure}
The closure of $A\sset\pw X$ in $\two_S^X$ is
\[\{Y\sset X:\forall Y_0\sset Y\text{ finite }\exists Z\in A\,Y_0\sset Z\}.\]
\end{Rem}

One more notational clarification: in \ref{item:5} below and elsewhere, we view $k\in\omega$ as a von~Neumann numeral,
that is, $k=\{0,\dots,k-1\}$.
\begin{Def}\th\label{def:cc}
A set $C\sset\pmf$ is a \emph{pmf clone} if the following hold for all $n,m,r,n',m'\in\omega$:
\def\theenumi{(\Roman{enumi})}
\begin{enumerate}
\item\label{item:1}
$C\cap\pmf_{n,m}$ is topologically closed as a subset of~$\two_S^{B^n\times B^m}$.
\item\label{item:2}
$C$ contains the identity function $\id_n\colon B^n\to B^n$.
\item\label{item:3}
$C$ is closed under composition: if $f\colon B^n\TO B^m$ and $g\colon B^m\TO B^r$ are in~$C$, then so is the pmf
$g\circ f\colon B^n\TO B^r$ given by 
\[(g\circ f)(x)\approx z\iff\exists y\in B^m\,(f(x)\approx y\land g(y)\approx z).\]
\item\label{item:4}
If $f\colon B^n\TO B^m$ and $g\colon B^{n'}\TO B^{m'}$ are in~$C$, then so is $f\times g\colon B^{n+n'}\TO B^{m+m'}$,
where $(f\times g)(x,x')\approx\p{y,y'}$ iff $f(x)\approx y$ and $g(x')\approx y'$.
\end{enumerate}
A class $D\sset\wgt$ is a \emph{weight coclone} if it satisfies the following conditions for any $k,k'\in\omega$:
\def\theenumi{(\Alph{enumi})}
\begin{enumerate}
\item\label{item:5}
If $w\colon B^k\to\strc M$ is in~$D$, and $\roo\colon k\to k'$, the weight $w\circ\tilde\roo\colon B^{k'}\to\strc M$ is in~$D$,
where $\tilde\roo(x^0,\dots,x^{k'-1})=\p{x^{\roo(0)},\dots,x^{\roo(k-1)}}$.
\item\label{item:6}
If $w\colon B^k\to\strc M$ is in~$D$, and $\fii\colon\strc M\to\strc M'$ is a pomonoid homomorphism (not necessarily
onto), then the weight $\fii\circ w\colon B^k\to\strc M'$ is in~$D$.
\item\label{item:7}
If $w_\alpha\colon B^k\to\strc M_\alpha$ is in~$D$ for every $\alpha\in I$, the weight $w\colon B^k\to\prod_{\alpha\in
I}\strc M_\alpha$ defined by $w(x)=\p{w_\alpha(x):\alpha\in I}$ is in~$D$.
\item\label{item:8}
If $w\colon B^k\to\strc M$ is in~$D$, and $\strc M'\sset\strc M$ is a submonoid including the image of~$w$, then $w\colon B^k\to\strc M'$ is
in~$D$.
\end{enumerate}
\end{Def}

\begin{Rem}\th\label{rem:least}
The smallest pmf clone is the clone $C_{\min}$ consisting of all subidentity partial functions $f\colon B^n\TO B^n$,
$f\sset\id_n$. Its dual coclone $\inv(C_{\min})$ is $\wgt$.

The index set~$I$ in condition~\ref{item:7} may be empty, in which case the result is the weight $w\colon B^k\to\strc1$
into the trivial pomonoid. In particular, every weight coclone is nonempty. The smallest weight coclone~$D_{\min}$
consists of all \emph{trivial weights} $\cst_1\colon B^k\to\strc M$, i.e., constant weights mapping to the unit of~$\strc
M$. The dual of~$D_{\min}$ is $\pol(D_{\min})=\pmf$.
\end{Rem}
\begin{Rem}\th\label{rem:ualg}
Preservation of weights by pmf can be put in the framework of ordered universal algebra in the following way. For fixed
$k\in\omega$, let $\Sigma_k$ denote the signature of pomonoids expanded with a set of extra constants $\{c_u:u\in B^k\}$. A
weight $w\colon B^k\to\strc M$ can be represented by a $\Sigma_k$-structure, namely $\strc M$ expanded with $w(u)$ as a
realization of the constant~$c_u$ for each $u\in B^k$. Let us denote this structure as~$\p{\strc M,w}$. We
see from~\eqref{eq:1} that for any $C\sset\pmf$, there is a set $I_C$ of inequalities between closed (variable-free) $\Sigma_k$-terms
such that $C\pres w$ iff $\p{\strc M,w}\model I_C$. In this way, $\inv(C)\cap\wgt_k$ becomes a partially ordered variety.
\end{Rem}
\begin{Rem}\th\label{rem:orders}
Let $D$ be a weight coclone, $\strc M$ a monoid, and $w\colon B^k\to\strc M$. Then the set of partial orders~$\le$ such that
$w\colon B^k\to\p{\strc M,\le}$ is in~$D$, is either empty, or a principal filter in the poset of partial orders
compatible with~$\strc M$; that is, it is closed under order extensions (by~\ref{item:6}), and nonempty intersections (by
\ref{item:7} and~\ref{item:8}, using the diagonal embedding of $\p{\strc M,\bigcap_\alpha{\le}_\alpha}$ in $\prod_\alpha\p{\strc M,{\le}_\alpha}$).

This can be even more naturally stated for invariant preorders in place of partial orders: that is, for any
pomonoid~$\strc M$, and $w\colon B^k\to\strc M$, the set of invariant preorders~$\preceq$ such that $w\colon
B^k\to\strc M/{\preceq}$ is in~$D$, is a principal filter in $\ocon\strc M$. Indeed, it is just $\ocon_{D_k}\p{\strc
M,w}$, where $D_k$ denotes $D\cap\wgt_k$ as a $\Sigma_k$-povariety in the set-up of \th\ref{rem:ualg}.
\end{Rem}

The main results of this section state that the notions of pmf clones and weight coclones faithfully describe
the closure systems of our Galois connection.

\begin{Thm}\th\label{thm:main-pmf}
Galois-closed sets of pmf are exactly the pmf clones.
\end{Thm}
\begin{Pf}
First, we show that any Galois-closed set of pmf is a pmf clone. Assume that $\pol(D)$ is such a set.
Let $n,m,k\in\omega$, $f$ be in the topological closure of $\pol(D)\cap\pmf_{n,m}$, $w\in D\cap\wgt_k$, and $a\in
B^{k\times n}$, $b\in B^{k\times m}$  be
such that $f(a^j)\approx b^j$ for all $j<k$. There exists $f'\in\pol(D)$ in the basic open neighbourhood $\{f':\forall
j<k\,f'(a^j)\approx b^j\}$ of~$f$; since $f'\pres w$, we obtain
\[\prod_{i<n}w(a_i)\le\prod_{i<m}w(b_i).\]
This shows that $\pol(D)$ satisfies~\ref{item:1}. Reflexivity and transitivity of~$\le$ readily implies
\ref{item:2} and~\ref{item:3}. Finally, if $f,g\in\pol(D)$, $w\in D$, and $(f\times g)(a^j,a'^j)\approx\p{b^j,b'^j}$
for each $j<k$, then
\[\prod_{i<n+n'}w((a,a')_i)=\Bigl(\prod_{i<n}w(a_i)\Bigr)\Bigl(\prod_{i<n'}w(a'_i)\Bigr)
\le\Bigl(\prod_{i<m}w(b_i)\Bigr)\Bigl(\prod_{i<m'}w(b'_i)\Bigr)=\prod_{i<m+m'}w((b,b')_i)\]
as $\cdot$ is nondecreasing in both arguments, which verifies~\ref{item:4}.

On the other hand, let $C$ be a pmf clone; we will prove $C$ is Galois-closed. In fact, we will construct a canonical
sequence of weights $\{w_k:k<\omega\}$ that characterize~$C$. For any $k<\omega$, let $\strc F_k=\p{F_k,1,\cdot}$ be
the monoid freely generated by $B^k$ (i.e., the monoid of finite words over alphabet~$B^k$), and define a relation
$\lesssim$ on~$F_k$ by
\begin{equation}\label{eq:20}
a_0\dots a_{n-1}\lesssim b_0\dots b_{m-1}\iff\exists g\in C\cap\pmf_{n,m}\,\forall j<k\,g(a^j)\approx b^j.
\end{equation}
Conditions \ref{item:2} and~\ref{item:3} imply that $\lesssim$ is a preorder, and \ref{item:4} shows that
$x\lesssim y$ implies $xz\lesssim yz$ and $zx\lesssim zy$. Thus, the relation $x\sim y\EQ x\lesssim y\land
y\lesssim x$ is a congruence on~$\strc F_k$, and $\strc M_k=\p{M_k,1,\cdot,\le}:=\p{F_k,1,\cdot,\lesssim}/{\sim}$ is
a pomonoid. Let $w_k\colon B^k\to\strc M_k$ be the quotient map composed with the natural inclusion of $B^k$ in~$\strc
F_k$. For any $a\in B^{k\times n}$, we have
\[\prod_{i<n}w_k(a_i)=(a_0\dots a_{n-1})/{\sim},\]
thus
\begin{equation}\label{eq:19}
\prod_{i<n}w_k(a_i)\le\prod_{i<m}w_k(b_i)\iff a_0\dots a_{n-1}\lesssim b_0\dots b_{m-1}.
\end{equation}
This implies immediately $C\pres w_k$ for all $k<\omega$.

Now, assume $f\in\pmf_{n,m}$ is in $\pol(\inv(C))$, which implies $f\pres\{w_k:k<\omega\}$; we want to show $f\in C$.
By \ref{item:1} and \th\ref{rem:closure}, it suffices to verify that for any $k<\omega$, $a\in B^{k\times n}$, and
$b\in B^{k\times m}$, if $f(a^j)\approx b^j$ for all $j<k$, then there exists $g\in
C\cap\pmf_{n,m}$ such that $g(a^j)\approx b^j$ for all $j<k$. This is indeed true: if $f(a^j)\approx b^j$ for all $j<k$, then $f\pres w_k$ implies $a_0\dots a_{n-1}\lesssim b_0\dots b_{m-1}$
by~\eqref{eq:19}, hence the required $g$ exists by~\eqref{eq:20}.
\end{Pf}
\begin{Thm}\th\label{thm:main-wgt}
Galois-closed classes of weights are exactly the weight coclones.
\end{Thm}
\begin{Pf}
Let $\inv(C)$ be any Galois-closed class of weights, we will show it is a weight coclone.
In view of \th\ref{rem:ualg}, for each fixed~$k$, $\inv(C)\cap\wgt_k$ is a povariety, and as such it is closed under
homomorphic images, products, and substructures. Moreover, since it is axiomatized by variable-free inequalities on
top of the theory of pomonoids, it is also closed under embedding into structures that are expansions of pomonoids,
whence under non-surjective homomorphisms into such structures. This shows that $\inv(C)$ satisfies conditions  \ref{item:6}--\ref{item:8}.
As for~\ref{item:5}, let $f\in C\cap\pmf_{n,m}$, and
$a'\in B^{k'\times n}$, $b'\in B^{k'\times m}$ be such that $f(a'^j)\approx b'^j$. Define $a\in B^{k\times n}$ and
$b\in B^{k\times m}$ by putting $a_i^j=a'^{\roo(j)}_i$,
$b_i^j=b'^{\roo(j)}_i$ for $j<k$. Then $f\pres w$ implies
\[\prod_{i<n}w(\tilde\roo(a'_i))=\prod_{i<n}w(a_i)\le\prod_{i<m}w(b_i)=\prod_{i<m}w(\tilde\roo(b'_i)).\]

Conversely, we prove that any weight coclone~$D$ is Galois-closed. Let $v\colon B^k\to\strc N$ be in $\inv(\pol(D))$,
we will show $v\in D$. Let $I$ denote the set of
pairs $\alpha=\p{(a_i^j)_{i<n}^{j<k},(b_i^j)_{i<m}^{j<k}}$ such that
\[\prod_{i<n}v(a_i)\nleq\prod_{i<m}v(b_i).\]
For each such~$\alpha$, let $f_\alpha\colon B^n\TO B^m$ denote the pmf $\{\p{a^j,b^j}:j<k\}$. Since $f_\alpha\npres v\in\inv(\pol(D))$, there
is $w'_\alpha\colon B^{k'}\to\strc M_\alpha$ in~$D$ such that $f_\alpha\npres w'_\alpha$. By the definition of~$f_\alpha$,
this means there is $\roo\colon k'\to k$ such that
\[\prod_{i<n}w'_\alpha(a^{\roo(0)}_i,\dots,a^{\roo(k'-1)}_i)\nleq
  \prod_{i<m}w'_\alpha(b^{\roo(0)}_i,\dots,b^{\roo(k'-1)}_i).\]
Thus, $w_\alpha=w'_\alpha\circ\tilde\roo\colon B^k\to\strc M_\alpha$, which is in~$D$ by~\ref{item:5}, satisfies
\begin{equation}\label{eq:2}
\prod_{i<n}w_\alpha(a_i)\nleq\prod_{i<m}w_\alpha(b_i).
\end{equation}
Let $\strc M'=\prod_{\alpha\in I}\strc M_\alpha$, $w'\colon B^k\to\strc M'$ be as in~\ref{item:7}, $\strc M$ be the
submonoid of~$\strc M'$
generated by $\rng(w')$, and $w\colon B^k\to\strc M$ be $w'$ reconsidered as a mapping to~$\strc M$. We have $w\in D$ by
\ref{item:7} and~\ref{item:8}.

If $\prod_{i<n}w(a_i)$ and $\prod_{i<m}w(b_i)$ are two elements of~$M$ such that
\begin{equation}\label{eq:3}
\prod_{i<n}w(a_i)\le\prod_{i<m}w(b_i),
\end{equation}
we must have $\prod_{i<n}v(a_i)\le\prod_{i<m}v(b_i)$: otherwise
$\alpha=\p{(a_i^j),(b_i^j)}\in I$, thus \eqref{eq:2} contradicts~\eqref{eq:3}. It follows that
\[\fii\Bigl(\prod_{i<n}w(a_i)\Bigr)=\prod_{i<n}v(a_i)\]
is a well-defined pomonoid homomorphism $\fii\colon\strc M\to\strc N$, hence $v=\fii\circ w\in D$ by~\ref{item:6}.
\end{Pf}

For any set $C\sset\pmf$, the least Galois-closed set containing~$C$ is $\pol(\inv(C))$. By \th\ref{thm:main-pmf},
$\pol(\inv(C))$ is exactly the pmf clone generated by~$C$, that is, the closure of $C$ under conditions
\ref{item:1}--\ref{item:4}. This is, on the face of it, a rather opaque operation, as in principle we might need to
cycle through the individual closure conditions and iterate them over and over. In fact, we will see that it is enough
to close the set under each condition once, in a judiciously chosen order, and similarly for the Galois closure of
classes of weights.
\begin{Def}
If $C\sset\pmf$, let $\cls_\cup C$ denote the closure of~$C$ under directed unions, $\cls_\sset C$ the closure
under subfunctions, and $\cls_{\id}C$, $\cls_\circ C$, and $\cls_\times C$ the closure under \ref{item:2},
\ref{item:3}, and \ref{item:4}, respectively. If $D\sset\wgt$, $\cls_{\var}D$, $\cls_MD$, $\cls_PD$, and
$\cls_SD$ denote the closure of~$D$ under \ref{item:5}, \ref{item:6}, \ref{item:7}, and \ref{item:8}, respectively.
\end{Def}

The $M$ in $\cls_M$ stands for ``morphism''. While $\cls_M$, $\cls_S$, and $\cls_P$ are reminiscent of the $H$, $S$,
and $P$ closure operators from universal algebra, condition~\ref{item:6} also applies to non-surjective homomorphisms;
for this reason, we chose a different letter to forestall confusion.
\begin{Cor}\th\label{cor:closure}
$\pol(\inv(C))=\cls_\cup\cls_\sset\cls_\circ\cls_\times\cls_{\id}C$, and
$\inv(\pol(D))=\cls_M\cls_S\cls_P\cls_{\var}D$.
\end{Cor}
\begin{Pf}
The $\Sset$ inclusions are clear, and the proof of \th\ref{thm:main-wgt} shows directly that any $w\in\inv(\pol(D))$ is in
$\cls_M\cls_S\cls_P\cls_{\var}D$. Let $f\in\pol(\inv(C))$, and $C^+$ be the closure of $C$ under \ref{item:2},
\ref{item:3}, and~\ref{item:4}. The argument in \th\ref{thm:main-pmf} shows that any finite subfunction of~$f$ is included
in some $g\in C^+$, hence $\pol(\inv(C))\sset\cls_\cup\cls_\sset C^+$. Clearly, $\cls_\circ\cls_\times\cls_{\id}C$
contains $\cls_{\id}C$, and it is closed under composition. It is also closed under $\times$, as
\[(f_r\circ\dots\circ f_1)\times(g_s\circ\dots\circ g_1)=(\id_m\times g_s)\circ\dots\circ(\id_m\times g_1)
   \circ(f_r\times\id_n)\circ\dots\circ(f_1\times\id_n),\]
where $m$ is the arity of the output of~$f_r$, and $n$ of the input of~$g_1$. Thus, $C^+=\cls_\circ\cls_\times\cls_{\id}C$.
\end{Pf}

\section{Restricted cases}\label{sec:restricted-cases}

The generality of the Galois connection described in Section~\ref{sec:clones-coclones} reflects more what we \emph{can}
do than what is \emph{useful} to do. In potential applications, we may be interested in restricted classes of
pmf or weights, for example:
\begin{itemize}
\item We may want to discuss only bona fide functions $f\colon B^n\to B^m$ (total, univalued) rather than pmf.
\item In the context of reversible operations, we want them further to be bijective, and we require $n=m$. On the other
hand, in the classical case we require $m=1$.
\item In essentially any reasonable context, we want to allow permutation of variables.
\item Some readers may prefer to disallow pesky corner cases involving~$B^0$, and only deal with pmf $B^n\TO B^m$ and
weights $B^k\to\strc M$ where $n,m,k\ge1$.
\item We may need to impose extra closure conditions, such as closure under inverse, or under usage of ancillary inputs
(see below).
\end{itemize}
We would like to adapt our Galois connection to such restricted contexts. Some cases are very easy to handle as an
immediate consequence of our main theorem:

\begin{Exm}\th\label{exm:n-n}
Let us investigate the Galois connection induced by the preservation relation restricted to pmf $f\colon B^n\TO B^n$
(i.e., with the same number of inputs and outputs). The class of all such pmf itself forms a clone, say~$C_=$; its dual
is a coclone~$D_=$. It follows that the preservation relation restricted to~$f\in C_=$ induces a Galois connection
whose closed classes are exactly the pmf clones $C\sset C_=$ on the one side, and weight coclones $D\Sset D_=$ on the
other side. In order to complete the description, it only remains to determine~$D_=$, ideally by presenting a simple
generating set. This is given in \th\ref{prop:restr} below: $D_=$ is generated by the constant-$1$ weight
function $\cst_1\colon B^0\to\p{\N,0,+,{=}}$. (In fact, $D_=$ consists of all constant weight functions.) Thus, the closed classes of weights under this restricted Galois connection
are weight coclones that contain~$\cst_1$. 
\end{Exm}

Several similarly easy cases can be dealt with using \th\ref{prop:restr,prop:zero} below.

Other cases turn out to be more complicated. In particular, the very important restriction of the Galois connection to
\emph{total} functions requires substantial work, and we will tackle it in
Section~\ref{sec:totality-conditions}. Likewise, in Section~\ref{sec:ancillas} we investigate closure conditions
imposed on classes of permutations in the work of Aaronson, Grier, and Schaeffer~\cite{aa-gr-sch}, in particular
closure under ancillas.

For the next statement, recall that the Kronecker delta function is defined by
\[\delta(u,v)=\begin{cases}1&u=v,\\0&u\ne v.\end{cases}\]
For any pmf $f\colon B^n\TO B^m$, its inverse $f^{-1}\colon B^m\TO B^n$ is defined by $f^{-1}(x)\approx y\EQ
f(y)\approx x$. A pmf $f$ is injective (sometimes called left-unique) if $f(x)\approx y$ and $f(x')\approx y$ implies
$x=x'$.
\begin{Prop}\th\label{prop:restr}
Let $C=\pol(D)$ and $D=\inv(C)$. If $w\colon B^k\to\strc M$ is a weight, let $\strc M_{(w)}$ denote the submonoid of~$\strc
M$ generated by $\rng(w)$.
\begin{enumerate}
\item\label{item:9}
All $f\colon B^n\TO B^m$ in~$C$ satisfy $n\le m$ ($n\ge m$; $n=m$) iff $D$ contains the constant weight $\cst_1\colon
B^0\to\p{\N,0,+}$, where $\N$ is ordered by $\le$ ($\ge$; $=$; resp.).
\item\label{item:10}
All $f\in C$ are partial functions (injective; both) iff $D$ contains Kronecker $\delta\colon B^2\to\p{\two,1,\land}$,
where the monoid is ordered by $\le$ ($\ge$; $=$; resp.).
\item\label{item:11}
$C$ contains the swap function $B^2\to B^2$, $\p{x,y}\mapsto\p{y,x}$ (and consequently all variable permutations $B^n\to B^n$) iff $\strc M_{(w)}$ is commutative for every $w\in D$.
\item\label{item:21}
$C$ contains all (constant) functions $B^0\to B$ iff for every $w\colon B\to\strc M$ in $D$, $1$ is a bottom element in~$\strc M_{(w)}$.
\item\label{item:12}
$C$ contains the diagonal mappings $\Delta_m\colon B\to B^m$, $\Delta_m(x)=\p{x,\dots,x}$, for $m=0,2$ (and
consequently, for all $m\ge0$), and variable
permutations, iff for every $w\in D$, $\strc M_{(w)}$ is a meet-semilattice with a top element $\p{M_{(w)},\top,\land,\le}$.
\item\label{item:13}
$f\in C$ implies $f^{-1}\in C$ iff for every $w\colon B^k\to\p{\strc M,\le}$ in $D$, also $w\colon B^k\to\p{\strc M,=}$
is in~$D$.
\end{enumerate}
\end{Prop}
\begin{Pf}
\ref{item:9}--\ref{item:11} are straightforward.

\ref{item:21}: A straightforward computation shows that $C$ contains all functions $B^0\to B^1$ iff for every
$w\colon B^k\to\strc M$ in~$D$, and for every $u\in B$, we have $1\le w(u,\dots,u)$. This is equivalent to the special
case $k=1$ given in~\ref{item:21}, since for any $w\in D$ as above, the unary weight $w'\colon B^1\to\strc M$ defined by
$w'(u)=w(u,\dots,u)$ is also in~$D$ by condition~\ref{item:5}.

\ref{item:12}: As long as $\strc M_{(w)}$ is commutative, $w$ is
preserved by diagonal maps iff $\strc M_{(w)}$ satisfies $x\le x^m$ for all $m\ge0$; that is, $x\le1$ and $x\le x^2$. On
the one hand, this gives $xy\le x,y$. On the other hand, $z\le x,y$ implies $z\le z^2\le xy$. Thus, $\p{\strc M,\le}$ is a semilattice with meet~$\cdot$, and top element~$1$.

\ref{item:13}, right-to-left: if $f\in C$, and $w\in D$, we may assume $w\colon B^k\to\p{\strc M,=}$. Then
condition~\eqref{eq:1} is symmetric in $a$ and~$b$, hence $f\pres w$ implies $f^{-1}\pres w$. Left-to-right: let
$w\colon B^k\to\p{\strc M,\le}$ be in~$D$, and consider any $f\in C$ and $(a_i^j),(b_i^j)$ such that $f(a^j)\approx b^j$ for
all $j<k$. Since $f\pres w$, we have
\[\prod_{i<n}w(a_i)\le\prod_{i<m}w(b_i).\]
By assumption, also $f^{-1}\in C$, and $f^{-1}(b^j)\approx a^j$ for all $j<k$, hence
\[\prod_{i<m}w(b_i)\le\prod_{i<n}w(a_i).\]
Thus, $f$ preserves the weight $w\colon B^k\to\p{\strc M,=}$.
\end{Pf}

The various cases in \ref{item:9} and~\ref{item:10} of \th\ref{prop:restr} give rise to variants of the clone--coclone
duality where $\pmf$ is restricted to a smaller class in the same way as in \th\ref{exm:n-n}. Conversely,
\ref{item:11}--\ref{item:13} lead to variants where $\wgt$ is restricted, such as:
\begin{Exm}\th\label{exm:commut}
Consider the Galois connection induced by the preservation relation restricted to $\wgt_{\mathrm{comm}}$---weights
$w\colon B^k\to\strc M$ with $\strc M$ a commutative pomonoid. We claim that the closed classes of this Galois
connection are, on the one side, pmf clones containing the swap function, and on the other side, classes of weights $w\colon
B^k\to\strc M$ from $\wgt_{\mathrm{comm}}$ satisfying conditions \ref{item:5}--\ref{item:8}, where \ref{item:6} is
restricted to $\strc M'$ commutative (let us call them ``commutative coclones'').

Indeed, \th\ref{prop:restr}~\ref{item:11} shows immediately that $C\sset\pmf$ is of the form $\pol(D)$ for
some $D\sset\wgt_{\mathrm{comm}}$ iff it is a clone containing the swap. It is also clear that since
$\inv(C)$ is a coclone, $\inv(C)\cap\wgt_{\mathrm{comm}}$ is a commutative coclone. On the other hand, if $D$ is a
commutative coclone, then the class $D'$ of weights $w\colon B^k\to\strc M'$ such that $w\colon B^k\to\strc M$ is in
$D$ for some (commutative) subpomonoid $\strc M\sset\strc M'$ is a coclone, and $D=D'\cap\wgt_{\mathrm{comm}}$, thus
$D=\wgt_{\mathrm{comm}}\cap\inv(C)$ for some $C\sset\pmf$.
\end{Exm}

Let us also state explicitly one more case as it involves a nontrivial closure condition on~$C$:
\begin{Exm}\th\label{exm:unordered}
Let $\wgt_{\mathrm{unord}}$ be the class of functions $w\colon B^k\to\strc M$ where $\strc M$ is an (unordered) monoid,
identified with the trivially ordered pomonoid $\p{\strc M,{=}}$.
Consider the Galois connection induced by the preservation relation between $f\in\pmf$ and
$w\in\wgt_{\mathrm{unord}}$. Using \th\ref{prop:restr}~\ref{item:13}, and the argument in \th\ref{exm:commut}, it is
easy to see that the closed classes of this Galois connection are, on the one side, clones $C\sset\pmf$ closed under
taking inverses (for every $f\colon B^n\TO B^m$ in~$C$, also $f^{-1}\colon B^m\TO B^n$ is in~$C$), and on the other
side, classes $D\sset\wgt_{\mathrm{unord}}$ closed under \ref{item:5}--\ref{item:8}, where in~\ref{item:6},
$\fii\colon\strc M\to\strc M'$ is a monoid homomorphism.
\end{Exm}

The restrictions of the Galois connection corresponding to the cases of \th\ref{prop:restr} can be combined where it
makes sense.

Although quite similar in spirit to \th\ref{prop:restr}, we discuss the following restrictions separately.
\pagebreak[2]
\begin{Prop}\th\label{prop:zero}
Let $C=\pol(D)$ and $D=\inv(C)$. Let us say that a pmf $f\colon B^n\TO B^m$ has the \emph{right shape} if $n>0$ ($m>0$;
both).
\begin{enumerate}
\item\label{item:24}
$C$ is generated by a set of pmf of the right shape iff all $f\in C$ are of the right shape except for $\id_0$ and $\nul\colon
B^0\TO B^0$ iff $D$ contains the constant weight $\cst_0\colon B^0\to\p{\two,1,\land}$
with the monoid ordered by~$\le$ ($\ge$; $=$; resp.).
\item\label{item:25}
$D$ is generated by $w\colon B^k\to\strc M$ with $k>0$ iff all nontrivial $w\colon B^k\to\strc M$ in~$D$ have $k>0$ iff
$C$ contains the empty pmf $\nul\colon B^0\TO B^1$ and $\nul\colon B^1\TO B^0$ (hence $\nul\colon B^n\TO B^m$ for all
$n,m$).
\end{enumerate}
\end{Prop}
\begin{Pf}
Straightforward.
\end{Pf}
\begin{Cor}\th\label{cor:zero}
\ \begin{enumerate}
\item\label{item:29}
The preservation relation restricted to pmf of the right shape induces a Galois connection
whose closed classes are, on the one side, sets of pmf of the right shape satisfying the appropriate restrictions of
\ref{item:1}--\ref{item:4}, and on the other side, weight coclones that include $\cst_0\colon B^0\to\p{\two,1,\land}$
ordered by~$\le$ ($\ge$; $=$; resp.).
\item\label{item:30}
The preservation relation restricted to weights $w\colon B^k\to\strc M$ with $k>0$ induces a Galois connection
whose closed classes are, on the one side, pmf clones that include $\nul\colon B^n\TO B^m$ for all $n,m\ge0$, and on
the other side, classes of said weights satisfying the appropriate restrictions of \ref{item:5}--\ref{item:8}.
\item\label{item:31}
Let us consider the preservation relation simultaneously restricted to pmf of the right shape, and to weights
$w\colon B^k\to\strc M$ with $k>0$. The closed classes of the induced Galois connection are, on the one side, classes
of pmf of the right shape satisfying the appropriate restrictions of \ref{item:1}--\ref{item:4}, and including 
all $\nul\colon B^n\TO B^m$ of the right shape; on the other side, classes of said weights satisfying the appropriate
restrictions of \ref{item:5}--\ref{item:8}, and including the constant weight $\cst_0\colon B^1\to\p{\two,1,\land}$
ordered by~$\le$ ($\ge$; $=$; resp.).
\end{enumerate}
\end{Cor}
\begin{Pf}
\ref{item:29} and~\ref{item:30} follow from \th\ref{prop:zero} in a similar way as in
\th\ref{exm:n-n,exm:commut,exm:unordered}. However, let us prove \ref{item:31} in more detail as the two conditions
from \th\ref{prop:zero} somewhat contradict each other, the weight in \th\ref{prop:zero}~\ref{item:24} being nullary.

Clearly, a closed set of pmf must satisfy the restricted versions of \ref{item:1}--\ref{item:4}, and
contain all empty pmf of the right shape by \th\ref{prop:zero}. Likewise, a closed class of weights satisfies the
restricted versions of \ref{item:5}--\ref{item:8}, and contains all weights $w\colon B^k\to\strc M$ with $k>0$ that are
in the coclone generated by $\cst_0\colon B^0\to\p{\two,1,\land}$ (ordered in the indicated fashion), one of which
being $\cst_0\colon B^1\to\p{\two,1,\land}$.

On the other hand, let $C$ be a set of pmf satisfying the conditions in the statement of~\ref{item:31}, and let $C'$ be
$C$ together with $\id_0\colon B^0\TO B^0$, and all empty pmf. Then $C'$ is a pmf clone, and by \th\ref{prop:zero}, its
dual $\inv(C')$ is generated by a class $D$ of weights $w\colon B^k\to\strc M$ with $k>0$. Thus, $C'=\pol(D)$, and $C$
is the restriction of $\pol(D)$ to the set of pmf of the right shape, as we added into $C'$ only pmf of wrong shapes.

Similarly, let $D$ be a class of weights satisfying the conditions in the statement of~\ref{item:31}. Let $D'$ denote
$D$ together with all weights $w\colon B^0\to\strc M$ such that the unique element $a\in\rng(w)$ is an idempotent, and
satisfies $a\le1$ ($1\le a$; nothing; respectively). Then $D'$ is a weight coclone: in particular, for any weight
$w\colon B^0\to\strc M$ as just described, the corresponding constant weight $B^1\to\strc M$ is already in~$D$ by
property~\ref{item:6}, as it factors through $\cst_0\colon B^1\to\p{\two,1,\land}$. By \th\ref{prop:zero},
$D'=\inv(C)$ for a set $C$ of pmf of the right shape. Again, the same holds for $D$ under the restricted preservation
relation as we only added nullary weights into~$D'$.
\end{Pf}

\begin{Rem}\th\label{rem:shape}
In general, nullary weights in $\inv(C)$ describe possible shapes of pmf in~$C$, where by
the shape of $f\colon B^n\TO B^m$ we mean the numbers $n,m$ of inputs and outputs. To see this, consider a weight
$w\colon B^0\to\strc M$. Let $a\in M$ be the single value of~$w$, and $\fii\colon\p{\N,0,+}\to\strc M$ the unique monoid
homomorphism mapping $1$ to~$a$. The order kernel $\oker(\fii)$ is an invariant preorder $\preceq$ on~$\N$, which faithfully
represents the relevant part of~$w$ in that the restricted-image weight $w\colon B^0\to\strc M_{(w)}$ (with the same polymorphisms) is isomorphic to the
weight $\cst_1\colon B^0\to\N/{\preceq}$ with single value~$1/{\preceq}$. Now, for any pmf $f\colon B^n\TO B^m$, we have
\[f\pres w\iff n\preceq m.\]

On a related note, recall from \th\ref{rem:orders} that for any pmf clone $C$, there is the smallest invariant preorder
$\preceq$ on~$\N$ such that $\cst_1\colon B^0\to\N/{\preceq}$ is in $\inv(C)$. It follows from the discussion above
that
\[n\preceq m\iff C\cap\pmf_{n,m}\ne\nul.\]
\end{Rem}

\begin{Rem}\th\label{rem:unary}
Case \ref{item:12} of \th\ref{prop:restr} directly gives a restricted version of the Galois connection with weights in
semilattices in the spirit of \th\ref{exm:commut}. However, the real reason we mention it is that it
can be used to recover the classical clone--coclone Galois connection, or more precisely, its version for partial
multifunctions $B^n\to B$ as given in the original paper by Geiger~\cite{geig} (we have yet to handle total
functions). We will describe the reduction now.

If a pmf clone~$C$ contains variable permutations and the diagonal maps $\Delta_m$, it also contains all projections
$\pi_{n,i}\colon B^n\to B$, $\pi_{n,i}(x_0,\dots,x_{n-1})=x_i$ (using $\Delta_0$). We claim that for any
$f\in\pmf_{n,m}$,
\begin{equation}\label{eq:4}
f\in C\iff\forall i<m\,\pi_{m,i}\circ f\in C.
\end{equation}
The left-to-right implication
follows from~\ref{item:3}, as $\pi_{m,i}\in C$. For the right-to-left implication, if $f_i=\pi_{m,i}\circ f\in C$ for
each $i<m$, then $f_0\times\dots\times f_{m-1}\colon B^{nm}\TO B^m$ is in~$C$. By using $\Delta_n$ and variable
permutations, $C$ includes the function $B^n\to B^{nm}$, $x\mapsto\p{x,\dots,x}$, hence also the pmf
$f'=\p{f_0,\dots,f_{m-1}}\colon B^n\TO B^m$ defined by
\[f'(x)\approx\p{y_0,\dots,y_{m-1}}\iff\forall i<m\,f_i(x)\approx y_i.\]
If $f$ were a function, then simply $f'=f$; for a general pmf~$f$, we still have that $f(x)\approx y$ implies
$f_i(x)\approx y_i$, hence $f\sset f'$, and $f\in C$ by~\ref{item:1}.

Thus, $C$ is determined by its unary-output fragment $C_1=\bigcupd_n(C\cap\pmf_{n,1})$. Now, $C_1$ satisfies \ref{item:1},
contains all $\pi_{n,i}$, and it is closed under composition in the sense that whenever it contains $g\colon B^{n'}\TO B$, and
$f_i\colon B^n\TO B$ for $i<n'$, it also contains $g\circ(f_0,\dots,f_{n'-1})\colon B^n\TO B$. Conversely, if
$C_1$ satisfies these three conditions, then \eqref{eq:4} defines a pmf clone~$C$ containing $\Delta_m$ and variable
permutations whose unary-output fragment is~$C_1$. Let us call such $C_1$ \emph{unary clones}.

On the dual side, if $w\colon B^k\to\strc M$ is a weight such that $\strc M_{(w)}$ is a meet semilattice with a top, we can
write $\strc M_{(w)}$ as a subdirect product of subdirectly irreducible such semilattices, i.e., $\p{\two,\top,\land,\le}$.
Such weights $w\colon B^k\to\two$ can be identified with relations~$r\sset B^k$, and it is easy to see that for
$f\colon B^n\TO B$, $f\pres r$ coincides with the classical preservation relation as in
Section~\ref{sec:examples}. Weight coclones~$D$ corresponding to unary clones are thus determined by classes $D_1$ of
relations $r\sset B^k$. The class~$D_1$ is closed under variable manipulations as in~\ref{item:5}. Conditions
\ref{item:6}--\ref{item:8} boil down to the following: If $r_\alpha\sset B^k$ is in~$D_1$ for all $\alpha\in I$, and
$F$ is a filter on $\pw I$, the relation $r\sset B^k$ defined by
\begin{equation}\label{eq:17}
r(x)\iff\{\alpha\in I:r_\alpha(x)\}\in F
\end{equation}
is in~$D_1$ (this describes when the weight corresponding to~$r$ can be obtained as a homomorphic image of a
subsemilattice of a direct product of weights corresponding to~$r_\alpha$). Notice that principal filters~$F$ give
closure of~$D_1$ under intersections. In general, one can check that $D_1$ is closed under~\eqref{eq:17} iff it is
closed in the Hausdorff topology, and under (finite, hence infinite) intersections. We thus obtain a Galois connection
whose closed classes are unary clones on one side, and classes of relations~$D_1$ closed under~\ref{item:5} and
intersections, and topologically closed, on the other side. This recovers the results of~\cite{geig}.
\end{Rem}

\subsection{Totality conditions}\label{sec:totality-conditions}
In contrast to the set-up of partial multifunctions that we worked with so far, our two original motivating examples (classical clones of operations on a
set, and reversible operations) deal exclusively with \emph{total} functions, it is thus imperative to investigate how
the duality is affected if we impose this restriction. That is, the preservation relation $f\pres w$ induces a Galois
connection between classes of \emph{total} multifunctions, and classes of weights: can we determine what are the closed
classes in this connection? Notice that unlike conditions like injectivity that we handled in
\th\ref{prop:restr,prop:zero}, the class of all total pmf is not a clone in our sense (it is obviously not closed downwards),
which considerably complicates the answer.

Closed classes of total mf are easy to guess; it is straightforward to formulate a version of the notion of
pmf clones for total mf (we will do this in \th\ref{def:total}). Moreover, it will mesh fairly well with the theory of
pmf clones, in that a clone of total mf can be identified with the pmf clone it generates by closing it downwards.
The pmf clones~$C$ we get in this way are those with the property that each pmf in~$C$ extends to a total mf in~$C$.

Description of the corresponding coclones is more difficult. Recall that in the classical case, coclones on a finite
set are sets of relations closed under positive primitive definitions, and in particular, totality corresponds to
closure under \emph{existential quantification}. (This correspondence does not work that well on infinite---especially
uncountable---sets, and we will encounter similar cardinality difficulties, too.) A natural generalization of
existential quantification to weight functions is as follows: if $w\colon B^{k+1}\to\strc M$ is a weight
function, let us define a weight $w^+\colon B^k\to\strc M$ by
\begin{equation}\label{eq:9}
w^+(x^0,\dots,x^{k-1})=\sum_{u\in B}w(x^0,\dots,x^{k-1},u).
\end{equation}
However, first we need to make sense of the sum in~\eqref{eq:9}. If $B$ is finite, it is enough to stipulate that
$\strc M$ be a \emph{semiring}: a structure $\p{M,1,{\cdot},0,{+}}$ where $\p{M,1,{\cdot}}$ is a monoid, $\p{M,0,{+}}$ is a
commutative monoid, and the finite distributive laws
\begin{align*}
(x+y)z&=xz+yz&z(x+y)&=zx+zy\\
0z&=0&z0&=0
\end{align*}
hold. If $B$ is infinite, we will also need some sort of completeness to make sense of infinite sums. We introduce the
relevant notions below.

The above discussion of totality applies symmetrically to classes of \emph{surjective (onto)} pmf, and we will treat
both cases, as well as their combination, in parallel.

\begin{Def}\th\label{def:semiring}
A \emph{partially ordered semiring (posemiring)} is a structure $\p{M,1,{\cdot},0,{+},{\le}}$ such that
$\p{M,1,{\cdot},0,{+}}$ is a semiring, and $\p{M,1,{\cdot},{\le}}$ and $\p{M,0,{+},{\le}}$ are pomonoids.\footnote{We
warn the reader that while this terminology is convenient for our purposes, and fits well in the general framework of
partially ordered varieties, it clashes with another commonly used
definition whereby ordered semirings, and in particular rings, satisfy the implication ${x\le y}\to{xz\le yz}\land{zx\le
zy}$ only for $z\ge0$, so the multiplicative monoid is not a pomonoid. On the other hand, they are additionally
required to satisfy $0\le1$.}

When $\strc M$ is a posemiring, we will abuse the language to speak of weight functions $w\colon B^k\to\strc M$ with the
understanding that this refers to the multiplicative pomonoid $\p{M,1,{\cdot},{\le}}$.

A \emph{positive (negative) semiring} is a posemiring $\p{M,1,{\cdot},0,{+},{\le}}$ such that $0\le1$ ($1\le0$, resp.).
This in fact implies that $0$ is a smallest (largest, resp.)\ element of~$\strc M$.

A semiring is \emph{idempotent} if it satisfies $x+x=x$. The additive structure of an idempotent semiring is a semilattice,
thus it can be interpreted as a pomonoid in two ways: a \emph{$\lor$-semiring} is an idempotent semiring ordered so
that $+$ is~$\lor$, and a \emph{$\land$-semiring} is an idempotent semiring ordered so that
$+$ is~$\land$. Equivalently, a $\lor$-semiring ($\land$-semiring) is an idempotent positive (negative, resp.)\
semiring. An idempotent semiring is \emph{continuous} if either of these two partial orders
makes it a complete lattice, satisfying the infinite distributive laws
\begin{equation}\label{eq:10}
\Bigl(\sum_{\alpha\in I}x_\alpha\Bigr)z=\sum_{\alpha\in I}x_\alpha z,\qquad z\sum_{\alpha\in I}x_\alpha=\sum_{\alpha\in I}zx_\alpha.
\end{equation}
More generally, an idempotent semiring is \emph{$\kappa$-continuous} for a cardinal~$\kappa$, if $\sum_{\alpha\in
I}x_\alpha$ exists and satisfies~\eqref{eq:10} whenever $\lh I<\kappa$. Notice that every idempotent semiring is
$\omega$-continuous.

Continuous $\lor$-semirings are commonly known in the literature as \emph{unital quantales}.
\end{Def}

The next theorem is the main technical result in this section: it characterizes pmf clones~$C$ whose dual coclones are
closed under the $w^+$ operation in continuous $\lor$-semirings. It turns out that on the fundamental level, this
operation does not precisely correspond to totality, but to the possibility of extending the domains of finite pmf
in~$C$ by one new element. (This also agrees with the effects of existential quantification in the classical duality.)
Only if the base set~$B$ is additionally \emph{countable} (including finite), we can repeat this extension countably many times to obtain that
each finite pmf from~$C$ extends to a total mf in~$C$.
\begin{Thm}\th\label{prop:ext}
Let $C=\pol(D)$ and $D=\inv(C)$. The following are equivalent.
\begin{enumerate}
\item\label{item:17}
For every finite $f\colon B^n\TO B^m$ in~$C$, and $a\in B^n$, there exists $b\in B^m$ such that $f\cup\{\p{a,b}\}\in
C$.
\item\label{item:18}
Whenever $w\colon B^{k+1}\to\strc M$ is in~$D$, where $\strc M$ is a continuous $\lor$-semiring, then also $w^+\colon
B^k\to\strc M$ is in~$D$.
\end{enumerate}
In condition~\ref{item:18}, it would be enough to demand that $\strc M$ be $\lh B^+$-continuous.

The symmetric condition $\forall f\in C\,\forall b\,\exists a\,f\cup\{\p{a,b}\}\in C$ is similarly characterized using continuous $\land$-semirings.
\end{Thm}
\begin{Pf}

\ref{item:17}\txto\ref{item:18}:
Assume $w\in D$, we will show $C\pres w^+$. Let $f\in C$, $f\colon B^n\TO B^m$, and $(a_i^j)_{i<n}^{j<k}$, $(b_i^j)_{i<m}^{j<k}$ be such that $f(a^j)\approx b^j$ for all
$j<k$. We may assume that $f$ is finite. Fix $a^k\in B^n$, and let $b^k\in B^m$ be such that
$g=f\cup\{\p{a^k,b^k}\}\in C$. Since $g\pres w$,
\[\prod_{i<n}w(a_i^0,\dots,a_i^k)\le\prod_{i<m}w(b_i^0,\dots,b_i^k)\le\prod_{i<m}w^+(b_i^0,\dots,b_i^{k-1}).\]
As $a^k$ was arbitrary, distributivity gives
\begin{align*}
\prod_{i<n}w^+(a_i^0,\dots,a_i^{k-1})&=\prod_{i<n}\LOR_{u\in B}w(a_i^0,\dots,a_i^{k-1},u)\\
&=\LOR_{a^k\in B^n}\prod_{i<n}w(a_i^0,\dots,a_i^k)\\
&\le\prod_{i<m}w^+(b_i^0,\dots,b_i^{k-1}).
\end{align*}

\ref{item:18}\txto\ref{item:17}:
Assume for contradiction that $(a_i^j)_{i<n}^{j\le k}$, $(b_i^j)_{i<m}^{j<k}$ are such that the pmf
$f=\{\p{a^j,b^j}:j<k\}$ is in~$C$, but
$f\cup\{\p{a^k,b^k}\}\notin C$ for all $b^k\in B^m$. For each such $b^k$, we can pick $w_{b^k}\colon B^{k+1}\to
\strc M_{b^k}$ in~$D$ such that
\[\prod_{i<n}w_{b^k}(a_i^0,\dots,a_i^k)\nleq\prod_{i<m}w_{b^k}(b_i^0,\dots,b_i^k)\]
using~\ref{item:5}. Let $\strc M=\prod_{b^k\in B^m}\strc M_{b^k}$, and $w\colon B^{k+1}\to\strc M$ be as
in~\ref{item:7}, so that $w\in D$, and
\begin{equation}\label{eq:5}
\prod_{i<n}w(a_i^0,\dots,a_i^k)\nleq\prod_{i<m}w(b_i^0,\dots,b_i^k)
\end{equation}
for every $b^k\in B^m$.

Let $\ob M$ be the complete lattice of down-sets of~$\strc M$, which we make into a pomonoid $\p{\ob M,1\down,{\cdot},{\sset}}$
by putting $X\cdot Y=\{xy:x\in X,y\in Y\}\down$. It is easy to check that
\[Y\cdot\bigcup_{\alpha\in I}X_\alpha=\bigcup_{\alpha\in I}Y\cdot X_\alpha,\]
and similarly for multiplication from the right, hence $\ob{\strc M}=\p{\ob M,1\down,{\cdot},\nul,{\cup},{\sset}}$ is in fact a continuous $\lor$-semiring. The mapping $x\mapsto
x\down$ is an embedding of $\strc M$ into~$\ob{\strc M}$, whose composition with~$w$ is a weight $\ob w\colon
B^{k+1}\to\ob{\strc M}$ in~$D$ by~\ref{item:6}. By assumption, the weight $\ob w^+\colon
B^k\to\ob{\strc M}$ given by
\[\ob w^+(x^0,\dots,x^{k-1})=\bigcup_{u\in B}\ob w(x^0,\dots,x^{k-1},u)
 =\{w(x^0,\dots,x^{k-1},u):u\in B\}\down\]
is also in~$D$. We have
\[\prod_{i<n}w(a_i^0,\dots,a_i^k)\in\prod_{i<n}\ob w^+(a_i^0,\dots,a_i^{k-1}),\]
but \eqref{eq:5} implies
\[\prod_{i<n}w(a_i^0,\dots,a_i^k)\notin\Bigl\{\prod_{i<m}w(b_i^0,\dots,b_i^k):b^k\in B^m\Bigr\}\Big\downarrow
 =\prod_{i<m}\ob w^+(b_i^0,\dots,b_i^{k-1}),\]
thus
\[\prod_{i<n}\ob w^+(a_i^0,\dots,a_i^{k-1})\nsset\prod_{i<m}\ob w^+(b_i^0,\dots,b_i^{k-1}).\]
This contradicts $f\in C$.
\end{Pf}
\begin{Def}\th\label{def:total}
A \emph{total (surjective; total surjective) clone} is a set $C$ of total multifunctions (surjective pmf; total
surjective mf; resp.)\ that satisfies
\ref{item:1}--\ref{item:4}, where condition~\ref{item:1} is understood relative to the subspace of $\two_S^{B^n\times
B^m}$ of all total (surjective; both; resp.)\ pmf.

A \emph{total (surjective; total surjective) coclone} is a weight coclone that satisfies condition~\ref{item:18} of
\th\ref{prop:ext} (the dual condition for~$\land$; both; resp.).
\end{Def}
\pagebreak[2]
\begin{Lem}\th\label{lem:total}
Assume that $B$ is countable, and let $C$ be a pmf clone.
\begin{enumerate}
\item\label{item:32}
If $C$ satisfies condition~\ref{item:17} from \th\ref{prop:ext}, then for every finite $f\colon B^n\TO B^m$ in~$C$,
there exists a total mf $\ob f\colon B^n\TO B^m$ in~$C$ such that $f\sset\ob f$.
\item\label{item:33}
If $C$ satisfies condition~\ref{item:17} from \th\ref{prop:ext} and its symmetric condition, then for every finite
$f\colon B^n\TO B^m$ in~$C$, there exists a total surjective mf $\ob f\colon B^n\TO B^m$ in~$C$ such that $f\sset\ob f$.
\end{enumerate}
\end{Lem}
\begin{Pf}

\ref{item:32}:
The statement is trivial if $B^n=\nul$; otherwise, let us fix a (not necessarily injective) enumeration
$B^n=\{a^j:j\in\omega\}$. By repeated application of condition~\ref{item:17}, we build an increasing chain of finite
pmf
\begin{equation}\label{eq:18}
f=f_0\sset f_1\sset f_2\sset f_3\sset\cdots
\end{equation}
such that $f_j\in C$, and $a^j\in\dom(f_{j+1})$. Then $\ob f=\bigcup_{j\in\omega}f_j\in C$ is a total extension of~$f$.

\ref{item:33} is similar: we enumerate $B^n=\{a^j:j\in\omega\}$ and $B^m=\{b^j:j\in\omega\}$, and we construct a
chain~\eqref{eq:18} such that $a^j\in\dom(f_{2j+1})$, and $b^j\in\rng(f_{2j+2})$.
\end{Pf}

\begin{Cor}\th\label{cor:total}
The preservation relation induces a Galois connection between sets of total (surjective; total surjective) pmf, and
classes of weights. In this connection, the Galois-closed sets of pmf are exactly the total (surjective; total
surjective; resp.)\ clones. All Galois-closed classes of weights are
total (surjective; total surjective; resp.)\ coclones, and if $B$ is countable, the converse also holds.
\end{Cor}
\begin{Pf}
We will discuss the total case, the other two cases are completely analogous.

Since the class of weights is unrestricted, the Galois-closed sets of pmf of this connection are exactly the
intersections of pmf clones with the set $\mathrm{Tmf}$ of all total mf by \th\ref{thm:main-pmf}. These are just
the total clones: on the one hand, if $C$ is a clone, then $C\cap\mathrm{Tmf}$ clearly satisfies the closure
conditions \ref{item:1}--\ref{item:4} restricted to~$\mathrm{Tmf}$; on the other hand, if $C$ is a total clone,
let $\ob C$ be the set of all pmf $g\colon B^n\TO B^m$ such that every finite $f\sset g$ is included in some $h\in C$.
Then using \th\ref{cor:closure}, we see that $\ob C$ is a pmf clone (that is, it is the clone generated by~$C$), while
the restricted condition~\ref{item:1} ensures that $\ob C\cap\mathrm{Tmf}=C$.

Galois-closed classes of weights of the restricted connection are thus classes of the form $D=\inv(C)$, where $C$ is a
total clone. Clearly, any such $D$ is a weight coclone. Moreover, since $D=\inv(\ob C)$, and any finite pmf in~$\ob C$
is included in some $f\in C$, which is total, we see that the pmf clone~$\ob C$ satisfies condition~\ref{item:17} of
\th\ref{prop:ext}, thus $D$ is a total coclone.

Conversely, if $D$ is a total coclone, then $\pol(D)$ is a pmf clone satisfying condition~\ref{item:17} of
\th\ref{prop:ext}. Thus, if $B$ is countable, then $\pol(D)$ is generated by $C=\pol(D)\cap\mathrm{Tmf}$ by
\th\ref{lem:total}, and as such $D=\inv(C)$ is a closed class of the restricted Galois connection.
\end{Pf}

We remark that the Galois connections from \th\ref{cor:total} can be combined with the restrictions from
\th\ref{prop:restr,prop:zero}: for example, we obtain (for countable~$B$) a duality between total clones consisting of
total \emph{functions}, and total coclones that include the Kronecker delta weight $\delta\colon
B^2\to\p{\two,1,\land,{\le}}$.

\begin{Exm}\th\label{exm:unctbl}
Neither \th\ref{lem:total} nor \th\ref{cor:total} holds if we drop the requirement of $B$ being countable. (Another
version of the example also applies to \th\ref{thm:permut} below.)

Assume that $B$ is uncountable, and fix a dense linear order $<$ on~$B$ such that there are copies of~$\Q$ at both ends
of $\p{B,<}$. Let $C$ be the pmf clone generated by all strictly order-preserving partial functions $B\TO B$, the
diagonal functions $\Delta_m$, and variable permutations. (The last two are not really necessary; we only include them
so that the example can be realized as a clone of partial operations $B^n\TO B$ in the classical set-up of
Geiger~\cite{geig}, cf.\ \th\ref{rem:unary}.) Explicitly, $C$ consists of subfunctions of partial functions $f\colon
B^n\TO B^m$ of the form
\[f(x_0,\dots,x_{n-1})=\p{f_0(x_{i_0}),f_1(x_{i_1}),\dots,f_{m-1}(x_{i_{m-1}})},\]
where $i_0,\dots,i_{m-1}<n$, and each $f_j\colon B\TO B$ is a strictly order-preserving partial function.

The pmf clone~$C$ satisfies condition \ref{item:17} of \th\ref{prop:ext}: we can extend each $f_j$ separately if it is
finite, using the density of~$<$.

However, $C$ does not satisfy the conclusion of \th\ref{lem:total}. Using the properties of~$<$, we can find two
elements $a,b\in B$ such that $a\up$ is uncountable, and $b\up$ is countable. Then the partial function $f\colon B\TO
B$ mapping $a$ to~$b$ is in~$C$, but it has no total extension in~$C$, as a total strictly order-preserving
extension of $f$ would need to embed $a\up$ in~$b\up$, which is impossible.

Furthermore, let $D=\inv(C)$. By \th\ref{prop:ext}, $D$ is a total coclone. However, $D$ is not a Galois-closed class
in the Galois connection restricted to total mf from \th\ref{cor:total}: indeed, if it were, then $C=\pol(D)$ would be
generated as a pmf clone by a set of total mf (a total clone), and as such it would satisfy condition~\ref{item:32}
from \th\ref{lem:total} (cf.\ the proof of \th\ref{cor:total}). We have just seen this is not the case. (We remark
that $D$ is generated by relations $r\sset B^k$ definable without parameters in the structure $\p{B,<}$, identified with
weights in $\p{\two,1,\land,{\le}}$ as in \th\ref{rem:unary}.)
\end{Exm}

Because of the application to reversible operations, we are particularly interested in clones determined by
permutations. At least for finite~$B$, their invariants are easily seen to be closed under the $w^+$ operation even
when the target posemiring is not idempotent. We will now investigate this closure condition more closely.

\begin{Thm}\th\label{prop:match}
Assume $B$ is finite, and let $C=\pol(D)$ and $D=\inv(C)$. The following are equivalent.
\begin{enumerate}
\item\label{item:19}
For every $f\colon B^n\TO B^m$ in~$C$, there is an injective function $g\colon B^n\to B^m$ such that $f\cup g\in C$
(which implies $n\le m$ unless $\lh B\le1$).
\item\label{item:20}
For every $w\colon B^{k+1}\to\strc M$ in~$D$ where $\strc M$ is a positive semiring, the weight $w^+\colon B^k\to\strc
M$ is in~$D$.
\end{enumerate}
\end{Thm}
\begin{Pf}

\ref{item:19}\txto\ref{item:20}:
Assume $w\colon B^{k+1}\to\strc M$ is in~$D$, we will verify $C\pres w^+$. Let $f\colon B^n\TO B^m$ be in~$C$, and $(a_i^j)_{i<n}^{j<k}$, $(b_i^j)_{i<m}^{j<k}$ be such that $f(a^j)\approx b^j$.
We may assume that $f$ includes an injection $g\colon B^n\to B^m$. Since $f\pres w$, we have
\begin{align*}
\prod_{i<n}w^+(a_i)&=\sum_{u\in B^n}\prod_{i<n}w(a_i,u_i)\\
&\le\sum_{u\in B^n}\prod_{i<m}w(b_i,g(u)_i)\\
&=\sum_{v\in g[B^n]}\prod_{i<m}w(b_i,v_i)\\
&\le\sum_{v\in B^m}\prod_{i<m}w(b_i,v_i)\\
&=\prod_{i<m}w^+(b_i),
\end{align*}
using positivity of the semiring.

\ref{item:20}\txto\ref{item:19}:
Let $f\colon B^n\TO B^m$ in~$C$ be $\sset$-maximal; we need to show that $f$ includes an injection.
Using Hall's marriage theorem, it suffices to verify that $\lh{f[X]}\ge\lh X$ for every $X\sset B^n$. 

Let us fix an enumeration $f=\{\p{a^j,b^j}:j<k\}$. As in the proof of \th\ref{prop:ext}, we can find a weight $w\colon B^{k+1}\to\strc M$
in~$D$ such that
\begin{equation}\label{eq:6}
\prod_{i<n}w(a_i,c_i)\nleq\prod_{i<m}w(b_i,d_i)
\end{equation}
for all $\p{c,d}\in(B^n\times B^m)\bez f$. We consider the monoidal semiring $\N[\strc M]$, whose elements are formal sums
\[x=\sum_{u\in M}x_uu,\]
where $x_u\in\N$, and $x_u=0$ for all but finitely many $u\in M$. We define a relation $\le$ on~$\N[\strc M]$ by
\begin{equation}\label{eq:8}
x\le y\iff\forall U\sset M\text{ up-set}\colon\sum_{u\in U}x_u\le\sum_{u\in U}y_u.
\end{equation}
We readily see that $\le$ is a partial order, and $x\le y$ implies $x+z\le y+z$. If $x\le y$, $v\in M$, and $U\sset M$
is an up-set, then $U'=\{u:uv\in U\}$ is also an up-set, and we have
\[\sum_{u\in U}(xv)_u=\sum_{u'\in U'}x_{u'}\le\sum_{u'\in U'}y_{u'}=\sum_{u\in U}(yv)_u,\]
hence $xv\le yv$. The set $\{z\in\N[\strc M]:xz\le yz\}$ thus contains~$0$, $1$, and it is closed under $+$ and right
multiplication by elements of~$M$, hence it is all of~$\N[\strc M]$. Symmetrically, one can prove $zx\le zy$, thus
$\N[\strc M]$ is
a positive semiring. The natural inclusion $\strc M\sset\N[\strc M]$ is a pomonoid embedding, hence we can
treat $w$ as a weight $w\colon B^{k+1}\to\N[\strc M]$.

By assumption, the weight $w^+\colon B^k\to\N[\strc M]$ is also in~$D$, in particular
\begin{equation}\label{eq:7}
\sum_{c\in B^n}\prod_{i<n}w(a_i,c_i)=\prod_{i<n}w^+(a_i)\le\prod_{i<m}w^+(b_i)=\sum_{d\in B^m}\prod_{i<m}w(b_i,d_i).
\end{equation}
Let $X\sset B^n$, and
\[U=\Bigl\{\prod_{i<n}w(a_i,c_i):c\in X\Bigr\}\Big\uparrow.\]
Then \eqref{eq:6}, \eqref{eq:8} and~\eqref{eq:7} give
\[\lh X\le\biggl|\Bigl\{c\in B^n:\prod_{i<n}w(a_i,c_i)\in U\Bigr\}\biggr|
\le\biggl|\Bigl\{d\in B^m:\prod_{i<m}w(b_i,d_i)\in U\Bigr\}\biggr|
=\bigl|f[X]\bigr|.\]
Since $X$ was arbitrary, Hall's theorem implies there exists an injection $g\sset f$.
\end{Pf}

Symmetrically, \th\ref{prop:match} also has a version with injective surjective partial functions in place of
injective (total) functions, and negative semirings in place of positive semirings. More interestingly, we can combine
both:
\begin{Cor}\th\label{cor:match}
Assume $B$ is finite, and let $C=\pol(D)$ and $D=\inv(C)$. The following are equivalent.
\begin{enumerate}
\item\label{item:27}
For every $f\colon B^n\TO B^m$ in~$C$, there is a bijection $g\colon B^n\to B^m$ such that $f\cup g\in C$ (which
implies $n=m$ unless $\lh B\le1$).
\item\label{item:28}
For every $w\colon B^{k+1}\to\strc M$ in~$D$ where $\strc M$ is a posemiring, the weight $w^+\colon B^k\to\strc M$ is
in~$D$.
\end{enumerate}
\end{Cor}
\begin{Pf}

\ref{item:27}\txto\ref{item:28}:
When $g$ is onto, the calculation in the proof of \ref{item:19}\txto\ref{item:20} of \th\ref{prop:match} is sound even
if the posemiring~$\strc M$ is not positive.

\ref{item:28}\txto\ref{item:27}:
We obtain an injective function~$g$ by \th\ref{prop:match}.
Any coclone $D$ contains the constant-$1$ weight $\cst_1\colon B\to\p{\N,1,{\cdot},0,{+},{\ge}}$, where the target is a
(negative) posemiring. Thus by our assumption, $D$ also contains the weight
$\cst_1^+\colon B^0\to\N$ whose single value is~$\lh B$. Since $f\pres\cst_1^+$, we obtain $\lh B^n\ge\lh B^m$, hence any
injection $g\colon B^n\to B^m$ is onto.
\end{Pf}

Still having reversible operations in mind, we close Section~\ref{sec:totality-conditions} by formulating a natural
version of our Galois connection for classes of permutations $B^n\to B^n$: we want closed classes of permutations to be
groups for each fixed~$n$ (in particular, to be closed under inverse), and to always include variable permutations.
These are in fact the demands on closed classes of reversible operations imposed by Aaronson, Grier,
and Schaeffer~\cite{aa-gr-sch}, except for closure under the ancilla rule, which we will handle in the next section.
On the dual side, the constraints on permutation clones allow us to only consider weights in commutative unordered
monoids, which simplifies the set-up.

\begin{Def}\th\label{def:perm}
A set of permutations $C\sset\bigcupd_{n\in\omega}\sym(B^n)$ is a \emph{permutation clone} if every $C\cap\sym(B^n)$ is
a closed subgroup of~$\sym(B^n)$ under its natural Hausdorff topology, $C$ is closed under~$\times$, and contains all
variable permutations.

A \emph{permutation weight} is a weight $w\colon B^k\to\strc M$, where $\strc M$ is a commutative monoid, considered as
a trivially ordered pomonoid. A class~$D$ of permutation weights is a \emph{permutation coclone}, if it satisfies \ref{item:5}--\ref{item:8} (with~\ref{item:6}
restricted to $\strc M'$ commutative and trivially ordered), contains the weights $\cst_1\colon B^0\to\p{\N,0,+}$
and~$\delta\colon B^2\to\p{\two,1,\land}$ from
\th\ref{prop:restr}, and is closed under the following version of condition~\ref{item:18} of \th\ref{prop:ext}: if $w\colon B^{k+1}\to\strc M$ is
in~$D$, where $\strc M$ is a continuous idempotent commutative semiring, then $D$ also contains $w^+\colon B^k\to\strc M$. (If $B$ is
finite, we may state this condition more generally with arbitrary commutative semirings.)
\end{Def}

\begin{Thm}\th\label{thm:permut}
The preservation relation induces a Galois connection between sets of permutations, and classes of permutation weights.
In this connection, the Galois-closed sets of permutations are exactly the permutation clones. All Galois-closed
classes of permutation weights are permutation coclones, and if $B$ is countable, the converse also holds.
\end{Thm}
\begin{Pf}
If $D$ is a class of permutation weights, the set of all permutations in $\pol(D)$ is a permutation clone by
\th\ref{thm:main-pmf,prop:restr}.

Conversely, let $C$ be a permutation clone, and $h\in\sym(B^n)\bez C$. Let $\ob C$ be the set of pmf $g\colon
B^m\TO B^m$ such that every finite subset of~$g$ is contained in some $f\in C$. By \th\ref{cor:closure}, $\ob
C=\pol(\inv(C))$. We have $h\notin\ob C$, hence there exists a
weight $w\colon B^k\to\p{\strc M,{\le}}$ such that $\ob C\pres w$, and $h\npres w$. We may assume $\le$ is~$=$ by
condition~\ref{item:13} of \th\ref{prop:restr}, and that $\strc M$ is commutative by condition~\ref{item:11}, hence $w$ is in
fact a permutation weight.

If $C$ is a permutation clone, the class of permutation weights in $\inv(C)=\inv(\ob C)$ is a permutation coclone by
\th\ref{thm:main-wgt,prop:restr,prop:ext} (and \th\ref{prop:match} if we use the extended
definition for $B$ finite).

On the other hand, let $D$ be a permutation coclone, and $\ob
D=\inv(\pol(D))$. Using \th\ref{cor:closure}, $\ob D$ consists of weights $w\colon B^k\to\p{\strc M,\le}$ for which there
exists a (necessarily commutative) submonoid $\rng(w)\sset\strc M'\sset\strc M$ such that $w\colon B^k\to\p{\strc M',=}$ is in~$D$; in
particular, all permutation weights in~$\ob D$ are in~$D$. Thus, if $w\notin D$ is a permutation weight, there is
$f\in\pol(D)$ such that $f\npres w$; we may assume $f$ is finite. The description of~$\ob D$ and \th\ref{prop:restr}
implies that $\pol(D)$ is closed under $\cdot^{-1}$, and that $\ob D$ satisfies condition~\ref{item:18} of
\th\ref{prop:ext} (whence also
its symmetric version). Thus, if $B$ is countable, we can extend $f$ to a
total surjective $\ob f\in\pol(D)$ by \th\ref{lem:total}. By \th\ref{prop:restr}, $\ob f$ is also an injective
function, hence it is in fact a permutation, and still $\ob f\npres w$.
\end{Pf}

\subsection{Ancillas and masters}\label{sec:ancillas}

There is another natural closure condition on clones commonly employed in reversible computing that we have not
discussed yet: the use of \emph{ancilla inputs}. The idea is that when we want to compute a permutation $B^n\to
B^n$, we may, along with the given input elements $x_0,\dots,x_{n-1}$, also work with auxiliary elements (ancillas)
$x_n,\dots,x_{n+m-1}$ that are initialized to a fixed string $a\in B^m$, as long as we guarantee to return these
extra elements to their original value at the end of the computation.

More formally, the (total) ancilla rule allows to construct a permutation $g\colon B^n\to B^n$ from a permutation $f\colon
B^{n+m}\to B^{n+m}$ if there exists $a\in B^m$ such that for all $x\in B^n$,
\[f(x,a)=\p{g(x),a}.\]
The usefulness of this rule in reversible circuit classes stems from the facts that on the one hand, it is considered
available for implementation, and on the other hand, it makes construction of circuits much more flexible; in particular, there is no
other way of producing a reversible circuit with a smaller number of inputs than what we started with.

While it is suggested in \cite{aa-gr-sch} that ancillas are similar to fixing inputs to constants in classical circuit
classes, this is only a loose analogy; in our set-up, we can formulate both, and closure under the ancilla rule
turns out to behave rather differently from closure under substitution of constants. In fact, we already dealt with the
latter: due to closure under composition, a pmf clone is closed under fixing inputs to constants iff it contains all
constant functions $B^0\to B^1$, which is in our duality equivalent to a restriction on unary weights presented in
\th\ref{prop:restr}~\ref{item:21}.

The ancilla rule as such is somewhat difficult to fit in our framework. For one thing, the rule is
``semantic''\footnote{In the sense used in ``syntactic vs.\ semantic complexity classes''.} in that we need to know that
all $x\in B^n$ satisfy a certain property before being allowed to construct the new function. It is also not very clear
how the rule should be generalized outside permutations. In order to get started, we characterize below in
\th\ref{thm:anc} pmf clones that are closed under its modified form---a \emph{partial} ancilla rule which is always applicable (relying on no
semantic promises), at the expense that it may produce partial functions in a way which does not play nice with
totality conditions as in Section~\ref{sec:totality-conditions}.

\begin{Def}\th\label{def:canc}
An element $z$ of a pomonoid $\p{M,1,\cdot,\le}$ is \emph{right-order-cancellative} if $xz\le yz$ implies $x\le y$ for
every $x,y\in M$.
\end{Def}

\begin{Thm}\th\label{thm:anc}
Let $C=\pol(D)$ and $D=\inv(C)$. The following are equivalent.
\begin{enumerate}
\item\label{item:22}
For all $f\colon B^{n+1}\TO B^{m+1}$ in~$C$, and $c\in B$, the pmf $g\colon B^n\TO B^m$ defined by
\[g(x)\approx y \iff f(x,c)\approx\p{y,c}\]
is in~$C$.
\item\label{item:23}
For every $w\colon B^k\to\strc M$ in~$D$, there is $w'\colon B^k\to\strc M'$ in~$D$, and a homomorphism $\fii\colon\strc M'\to\strc M$ such
that $w=\fii\circ w'$, and the diagonal weights $w'(c,\dots,c)$ are right-order-cancellative in~$\strc M'$ for all $c\in B$.
\end{enumerate}
\end{Thm}
\begin{Pf}

\ref{item:23}\txto\ref{item:22}:
Let $f$ and~$g$ be as in~\ref{item:22}, and $w\colon B^k\to\strc M$ in~$D$. Let $w'\colon B^k\to\strc M'$ and~$\fii$ be
as in~\ref{item:23}; in particular, $w'(c^{(k)})$ is right-order-cancellative in~$\strc M'$, where
$c^{(k)}=\p{c,\dots,c}\in B^k$. If $(a_i^j)_{i<n}^{j<k}$, $(b_i^j)_{i<m}^{j<k}$ are such that $g(a^j)\approx b^j$, we have
\[\prod_{i<n}w'(a_i)\cdot w'(c^{(k)})\le\prod_{i<m}w'(b_i)\cdot w'(c^{(k)})\]
as $f\pres w'$, hence
\[\prod_{i<n}w'(a_i)\le\prod_{i<m}w'(b_i)\]
by cancellativity. Since $\fii$ is a pomonoid homomorphism and $w=\fii\circ w'$, this implies
\[\prod_{i<n}w(a_i)\le\prod_{i<m}w(b_i).\]

\ref{item:22}\txto\ref{item:23}:
Let $D'$ denote the class of weights $w'\colon B^k\to\strc M'$ in~$D$ such that $w'(c^{(k)})$ is
right-order-cancellative in~$\strc M'$ for all $c\in B$. Let $w_k\colon B^k\to\strc M_k$ be as in the proof of \th\ref{thm:main-pmf}. Condition~\ref{item:22} ensures that
$w_k\in D'$, and the proof of \th\ref{thm:main-pmf} shows that $C=\pol(\{w_k:k\in\omega\})$, hence
\[D=\inv(\pol(\{w_k:k\in\omega\}))=\cls_M\cls_S\cls_P\cls_{\var}\{w_k:k\in\omega\}\]
by \th\ref{cor:closure}. It is easy to see that $D'$ is closed under $\cls_S$, $\cls_P$, and~$\cls_{\var}$, hence
$D=\cls_M(D')$.
\end{Pf}

\emph{Reversible gate classes} as defined in~\cite{aa-gr-sch} are permutation clones closed under the \emph{total}
ancilla rule. We call them \emph{master clones} below, in accordance with the terminology used in~\cite{ej:revclass},
where a rudimentary form of our Galois connection was first presented.
We present in \th\ref{thm:ancilla} a convenient variant of our Galois connection for master clones on a finite base
set~$B$. As in \th\ref{def:perm}, we will only deal with unordered weights in commutative monoids.

If we take a master clone, and close it under subfunctions to generate a pmf clone in a minimal way, there is no reason
to expect the generated clone to be closed under the partial ancilla rule. However, we can get there in a roundabout
way: given a master clone~$C$, we close it under the partial ancilla rule (and subfunctions); we obtain a pmf clone
that is closed under the partial ancilla rule by definition, and crucially, that does not contain any new total
functions outside~$C$---this is precisely what the closure of~$C$ under the \emph{total} ancilla rule tells us.
(However, the pmf clone we get in this way does not have the property of extendability to total functions as in
\th\ref{prop:ext} or \th\ref{prop:match}, even though it was ``generated'' from a set of total functions!)

Thus, we will be able to describe master clones by weights satisfying a property that actually corresponds to the partial
ancilla rule as in \th\ref{thm:anc}. In the context of unordered commutative monoids, right order-cancellativity is
equivalent to plain cancellativity. We will in fact impose a stronger requirement, namely that the
monoid elements in question have an inverse; this makes the condition somewhat simpler, and as a technical advantage,
stable under homomorphisms. We rely here on the basic fact that cancellative elements of a commutative monoid can be
made invertible in a (commutative) extension of the monoid by a variant of the Grothendieck group construction:
\begin{Lem}\th\label{lem:groth}
If $\strc M$ is a commutative monoid, and $U\sset\strc M$ a set of cancellative elements, there exists a
commutative monoid $\strc N\Sset\strc M$ such that every $u\in U$ has an inverse in~$\strc N$.
\end{Lem}
\begin{Pf}
Let $\strc M_U$ be the submonoid of~$\strc M$ generated by~$U$. Since the elements of~$\strc M_U$ are cancellative, the
relation
\[\p{x,u}\sim\p{y,v}\iff xv=yu\]
on $\strc M\times\strc M_U$ is easily seen to be a congruence, thus we can form the quotient $\strc N=(\strc
M\times\strc M_U)/{\sim}$. The monoid $\strc M$ embeds in $\strc N$ via $x\mapsto\p{x,1}/{\sim}$, and for $u\in\strc
M_U$, $\p{u,1}/{\sim}$ has an inverse $\p{1,u}/{\sim}$.
\end{Pf}

Because of the complicated interference of the ancilla rules with totality conditions, we are not able to precisely
describe the closed classes of weights in the Galois connection for master clones; we only know some necessary
conditions.

\begin{Def}\th\label{def:ancilla}
Assume that $B$ is finite. A permutation clone~$C$ is a \emph{master clone} if it is closed under the total ancilla
rule: if $f\in C\cap\sym(B^{n+m})$, $a\in B^m$, and $g\in\sym(B^n)$ are such that $f(x,a)=\p{g(x),a}$ for all $x\in
B^n$, then $g\in C$.

A permutation weight $w\colon B^k\to\strc M$ is a \emph{master weight} if the diagonal weights $w(x^{(k)})$ are
invertible in~$\strc M$ for all $x\in B$. A \emph{master proto-coclone} is a class of master weights that contains
$\cst_1\colon B^0\to\p{\Z,0,+}$, $\delta\colon B^2\to\p{\two,1,\land}$, and is closed under
conditions \ref{item:5}--\ref{item:8} relative to the class of all master weights.
\end{Def}

\begin{Thm}\th\label{thm:ancilla}
If $B$ is finite, the preservation relation induces a Galois connection between sets of permutations and classes of
master weights such that Galois-closed sets of permutations are exactly the master clones. Galois-closed classes of
master weights are master proto-coclones.
\end{Thm}
\begin{Pf}
The set of permutations preserving a class of master weights is a permutation clone by \th\ref{thm:permut}, and it
satisfies the ancilla rule by the argument in the proof of \th\ref{thm:anc}. Likewise, the class of master weights
preserving a set of permutations is a master proto-coclone, being the intersection of a permutation
coclone with the class of all master weights.

Let $C$ be a master clone, and $h$ a permutation not in~$C$. Let $\ob C$ be the set of pmf $g\colon B^n\TO B^n$
such that there are $m\in\omega$, $a\in B^m$, and $f\in C\cap\sym(B^{n+m})$ such that
$f(x,a)=\p{y,a}$ whenever $g(x)\approx y$. We have $h\notin\ob C$. We can check easily that $\ob C$ is a pmf
clone closed under $\cdot^{-1}$, hence there is a permutation weight $w\colon B^k\to\strc M$ in $\inv(\ob C)$ such that
$h\npres w$ as in the proof of \th\ref{thm:permut}. Moreover, $\ob C$ satisfies condition~\ref{item:22} of
\th\ref{thm:anc}, hence we can assume that $w(x^{(k)})\in M$ is cancellative for every $x\in B$. By \th\ref{lem:groth},
we can embed $\strc M$ in a commutative monoid~$\strc N$ where every $w(x^{(k)})$ is invertible, thus $w\colon
B^k\to\strc N$ is a master weight.
\end{Pf}
\begin{Rem}
It is not clear what other properties do Galois-closed classes of master weights satisfy. Extending the argument in
the proof of \th\ref{thm:ancilla}, we can show that the following are equivalent for a master proto-coclone~$D$:
\begin{enumerate}
\item\label{item:34} $D$ is Galois-closed in the connection from \th\ref{thm:ancilla}.
\item\label{item:35} All master weights in the least permutation coclone containing~$D$ are already in~$D$.
\item\label{item:36} For every pmf $g\colon B^n\TO B^n$ such that $g\pres D$, there are $m\in\omega$, $f\in\sym(B^{n+m})$, and $a\in
B^m$ such that $f\pres D$, and $f(x,a)=\p{y,a}$ whenever $g(x)\approx y$.
\end{enumerate}
Unfortunately, these conditions are fairly opaque. We can at least infer additional necessary closure conditions on~$D$
from~\ref{item:35}, in particular the following: if $w\colon B^{k+l}\to\strc M$ is in $D$, where $\strc M$ is a
semiring, let $w^{+l}\colon B^k\to\strc M$ be the permutation weight defined by
\[w^{+l}(x)=\sum_{y\in B^l}w(x,y)\]
(i.e., the $w^+$~construction $l$~times iterated), let
\[u\sim v\iff\exists n\in\omega\,\exists x\in B^n\:u\prod_{i<n}w^{+l}(x_i^{(k)})=v\prod_{i<n}w^{+l}(x_i^{(k)})\]
be the least congruence on~$\strc M$ that makes the diagonal $w^{+l}$-weights cancellative, and let $\strc
N\Sset\strc M/{\sim}$ be the Grothendieck monoid from \th\ref{lem:groth} where they are made invertible. Then
$w^{+l}\colon B^k\to\strc N$ is in~$D$.

However, we see no evidence to suggest that this condition is sufficient.
\end{Rem}
\begin{Prob}\th\label{prob:master}
Describe Galois-closed classes of master weights by means of transparent closure conditions.
\end{Prob}

\section{Subdirectly irreducible weights}\label{sec:subd-irred-weights}

Coclones with weights in arbitrary pomonoids, as we discussed so far, are convenient for the abstract theory because of
their rich closure properties, but not so much for applications, as such invariants are absurdly large: any
nontrivial coclone is a proper class that includes pomonoids of arbitrary cardinality, most of which are clearly
redundant. In contrast, the classical clone--coclone duality only involves finitary relations on the base set~$B$,
which are small finite objects for finite~$B$.

In light of this, we would like to identify a smaller class of weights that suffice to characterize every clone.
The price we are willing to pay is that we resign on the idea of coclones having sensible closure conditions:
in particular, we drop closure under products, and in the process we will also lose closure properties such as in
Section~\ref{sec:totality-conditions}. In fact, we will in a sense
make closure conditions work backwards, to decompose any weight into small pieces that generate it.

Recall from \th\ref{rem:ualg} that any weight $w\colon B^k\to\strc M$ can be identified with an expansion $\p{\strc
M,w}$ of the pomonoid~$\strc M$ by extra constants describing the values of~$w$. If $D$ is a coclone, this makes
$D\cap\wgt_k$ into a povariety. Furthermore, we only care about the validity of variable-free inequalities in its
inequational theory.
This suggests two ways to obtain a small generating set in a coclone:
\begin{itemize}
\item We can restrict attention to weights $w\colon B^k\to\strc M$ such that $\strc M$ is generated by $\rng(w)$ as a monoid (that
is, $\p{\strc M,w}$ is $0$-generated). In particular, if $B$ is finite, this makes $\strc M$ finitely generated.
\item As in any povariety, $\p{\strc M,w}$ can be written as a subdirect product of subdirectly irreducible poalgebras; thus, we
can restrict attention to weights with the pomonoid~$\strc M$ subdirectly irreducible (subdirect irreducibility is unaffected by
expansions by constants).
\end{itemize}
Let us state the result explicitly for the record.
\begin{Prop}\th\label{prop:sdi-weights}
Every pmf clone~$C$ can be written in the form $C=\pol(D)$, where $D$ is a set of weights $w\colon B^k\to\strc M$ such
that the pomonoid $\strc M$ is subdirectly irreducible, and generated by~$\rng(w)$.
\noproof\end{Prop}

The same principle also applies to variants of the Galois connection with restricted classes of weights, as discussed
in Section~\ref{sec:restricted-cases}. Most cases of interest can be handled by the following generalization of
\th\ref{prop:sdi-weights}:
\begin{Prop}\th\label{prop:sdi-weights-gen}
Let $Q$ be a poquasivariety of pomonoids, and $P\sset\pmf$. Consider the Galois connection induced by the preservation
relation between pmf $f\in P$, and weights $w\colon B^k\to\strc M$ with $\strc M\in Q$. Then every Galois-closed set of
pmf is of the form $P\cap\pol(D)$, where $D$ is a set of weights as above with $\strc M$ subdirectly irreducible
relative to~$Q$, and generated by~$\rng(w)$.

Even more generally, we can use possibly different quasivarieties $Q_k$ for each $k\in\omega$.
\noproof\end{Prop}
The various restrictions on weights mentioned in
\ref{exm:n-n}--\ref{rem:shape}, as well as permutation weights from \th\ref{def:perm}, can be defined using
quasi-inequalities, hence they are in the scope of \th\ref{prop:sdi-weights-gen}. (Allowing $Q_k$ to depend on~$k$ is
useful e.g.\ if we want to disable nullary weights as in \th\ref{cor:zero}.) Notice that we have already implicitly
used a form of \th\ref{prop:sdi-weights-gen} in \th\ref{rem:unary}.

The one remaining exception is the class of master weights from \th\ref{def:ancilla}: these do not form a quasivariety,
as the invertibility of the diagonal weights $w(x,\dots,x)$ needs an existential quantifier to state. We can adapt our
approach to this case anyway: we can assume that
$\strc M$ is generated by $\rng(w)$ together with inverses of the diagonal weights, which again makes it finitely
generated if $B$ is finite; and we can restrict attention to $\strc M$ subdirectly irreducible (in the ordinary
algebraic sense).

These considerations are particularly useful for classes of weights whose finitely generated subdirectly irreducible
pomonoids are finite, in which case we can characterize all clones by finite invariants if $B$ is finite.
In particular, this happens for \emph{permutation weights}, which can describe all permutation clones, and more
generally, all pmf clones containing variable permutations and being closed under inverses.

The fact that finitely generated subdirectly irreducible commutative monoids (or semigroups) are finite was proved by
Mal'cev~\cite{malc:semigp}, and their structure has been investigated by Schein~\cite{schein} and
Grillet~\cite{grillet}. For completeness, we include a simplified description below.
\begin{Def}\th\label{def:sdisemigp}
A \emph{nilsemigroup} is a semigroup $\p{N,\cdot}$ with an absorbing element~$0$ (i.e., $0x=x0=0$ for all $x\in N$) such that for every $x\in N$, there is
$n\in\N_{>0}$ such that $x^n=0$.

Let $\p{G,1,\cdot}$ be an abelian group, $\p{\Omega,\cdot,0}$ a commutative nilsemigroup, and $\Omega^1$ the monoid
$\Omega\cup\{1\}$. A \emph{factor set} on~$\Omega$ with values in~$G$ is
$\sigma=\{\sigma_{\alpha,\beta}:\alpha,\beta\in\Omega^1,\alpha\beta\ne0\}\sset G$ satisfying
\begin{align*}
\sigma_{\alpha,\beta}&=\sigma_{\beta,\alpha},\\
\sigma_{\alpha,1}&=1,\\
\sigma_{\alpha,\beta}\,\sigma_{\alpha\beta,\gamma}&=\sigma_{\alpha,\beta\gamma}\,\sigma_{\beta,\gamma}
\end{align*}
for all $\alpha,\beta,\gamma\in\Omega^1$ such that $\alpha\beta\gamma\ne0$. We define a commutative monoid
$[\Omega,G,\sigma]$ whose underlying set is
\[\{0\}\cup\{\p{g,\alpha}:g\in G,\alpha\in\Omega^1,\alpha\ne0\},\]
with unit $\p{1,1}$, absorbing element~$0$, and multiplication of nonzero elements defined by
\[\p{g,\alpha}\p{h,\beta}=\begin{cases}
0&\alpha\beta=0,\\
\p{gh\sigma_{\alpha,\beta},\alpha\beta}&\text{otherwise.}
\end{cases}\]
Note that $[\Omega,G,\sigma]$ is a disjoint union of a copy of~$G$, and a nilsemigroup~$N$ which is an ideal of
$[\Omega,G,\sigma]$, such that the semigroup of orbits $N/{\equiv}$ is isomorphic to~$\Omega$, where $u\equiv v$ iff
$u=gv$ for some $g\in G$. In particular, $\Omega$ and~$G$ are determined by $[\Omega,G,\sigma]$ up to isomorphism. The
factor set~$\sigma$, however, is not: $[\Omega,G,\sigma]\simeq[\Omega,G,\sigma']$ if
\[\sigma'_{\alpha,\beta}=\sigma_{\alpha,\beta}u_\alpha u_\beta u_{\alpha\beta}^{-1}\]
for some $\{u_\alpha:\alpha\in\Omega\bez\{0\}\}\sset G$ with $u_1=1$.

$\Omega$ carries a canonical partial order, defined by $x\le y$ iff $x=uy$ for some $u\in\Omega^1$. Assume that $\Omega$ has
a unique minimal element~$\mu$. We say that $N$ is \emph{weakly irreducible} when the following
condition holds for all $\alpha,\beta\in\Omega$ (or $\Omega^1$): if
\[\{\tau\in\Omega:\tau\alpha=\mu\}=\{\tau\in\Omega:\tau\beta=\mu\},\]
and the mapping $\tau\mapsto\sigma_{\alpha,\tau}\sigma_{\beta,\tau}^{-1}$ is constant on $\{\tau\in\Omega:\tau\alpha=\mu\}$, then $\alpha=\beta$.
\end{Def}
\begin{Thm}[Grillet~\cite{grillet}]\th\label{thm:sdisemigp}
Let $G$ be a abelian group, $\Omega$ a finite commutative nilsemigroup, and $\sigma$ a factor set on~$\Omega$ with
values in~$G$ such that:
\begin{enumerate}
\item $G$ is trivial, or a cyclic group~$C_{p^k}$ of prime power order.
\item $\Omega$ is trivial, or it has a unique minimal element~$\mu$, and $N$ is weakly irreducible.
\end{enumerate}
Then $[\Omega,G,\sigma]$ is a finite subdirectly irreducible commutative monoid.

Conversely, every finite (= finitely generated) subdirectly irreducible commutative monoid is isomorphic to some
$[\Omega,G,\sigma]$ as above, or to $C_{p^k}$ (which equals $[0,C_{p^k},1]$ minus its zero element).
\end{Thm}
An example of a nilsemigroup with a unique minimal element is $\p{\{1,\dots,d\},\min\{d,x+y\}}$ for $d>1$:
here $d$ is a zero element, and $\mu=d-1$.

\section{Conclusion}\label{sec:conclusion}

We have demonstrated that the usual notion of clones of multiple-input single-output functions has a natural
generalization to multiple-output (partial multi-valued) functions. The clone--coclone Galois connection admits a
parallel generalization: multiple-output clones are determined by preservation of invariants in the shape of weighted
products in partially ordered monoids. Closed classes of invariants, generalizing the usual coclones, are likewise
characterized by suitable closure conditions resembling properties of varieties, and admit a form of decomposition into
subdirectly irreducible invariants.

We have also seen that the Galois connection is flexible enough to accommodate various natural modifications of the
set-up: for example, if we desire clones to be always closed under variable permutations, it suffices to switch to
commutative monoids; constraining clones to partial uni-valued functions corresponds to inclusion of a specific weight
as an invariant; and closure of clones under inverse can be obtained by using unordered monoids. Most importantly, we
can adapt the Galois connection to classes of total functions (multi-valued or uni-valued), in which case coclones get
closed under sums in certain semirings, generalizing the closure under existential quantification from the classical
case. By suitable restrictions on weights, we can also take care of the ancilla rule as employed in reversible
computing, although we could not exactly determine the corresponding closure conditions on the coclone side.

This work opens up the possibility to study multiple-output clones with tools of universal algebra, generalizing the
standard framework of clones to a set-up where domains and codomains of operations receive a fully symmetrical treatment.

\bibliographystyle{mybib}
\bibliography{mybib}
\end{document}
